\documentclass[12pt]{amsart}
\usepackage{latexsym, amssymb, amsfonts}
\newenvironment{pf}{\proof[\proofname]}{\endproof}
\theoremstyle{plain}
\newtheorem{Thm}{Theorem}[section]
\newtheorem{Cor}[Thm]{Corollary}
\newtheorem{Main}[Thm]{Main Theorem}

\newtheorem{Prop}[Thm]{Proposition}
\newtheorem{Lem}[Thm]{Lemma}
\newtheorem{Cl}[Thm]{Claim}

\theoremstyle{definition}
\newtheorem{Def}[Thm]{Definition}
\newtheorem{Rem}[Thm]{Remark}

\newtheorem{Emp}[Thm]{}
\newtheorem{Ex}[Thm]{Example}

\newtheorem{Not}[Thm]{Notation}

\numberwithin{equation}{section}

\newcommand{\B}[1]{\mathbb#1}
\newcommand{\cal}[1]{\mathcal{#1}}
\newcommand{\C}[1]{\cal#1}

\newcommand{\isom}{\overset {\thicksim}{\to}}

\newcommand{\hra}{\hookrightarrow}
\newcommand{\wt}{\widetilde}
\newcommand{\Gm}{\Gamma}

\newcommand{\bo}{\mathbf}

\newcommand{\gm}{\gamma}
\newcommand{\dt}{\delta}
\newcommand{\Dt}{\Delta}
\newcommand{\bs}{\backslash}

\newcommand{\m}{^{\times}}

\newcommand{\al}{\alpha}

\newcommand{\an}{^{an}}

\newcommand{\rs}[1]{Section \ref{S:#1}}
\newcommand{\rl}[1]{Lemma \ref{L:#1}}

\newcommand{\rp}[1]{Proposition \ref{P:#1}}
\newcommand{\rr}[1]{Remark \ref{R:#1}}
\newcommand{\rc}[1]{Construction \ref{C:#1}}
\newcommand{\re}[1]{\ref{E:#1}}
\newcommand{\rco}[1]{Corollary \ref{C:#1}}

\newcommand{\rt}[1] {Theorem \ref{T:#1}}

\newcommand{\sm}{\smallsetminus}
\newcommand{\be}{\infty}

\newcommand{\pr}{\text{pr}}

\newcommand{\Spec}{\text{Spec}}
\newcommand{\Aut}{\text{Aut}}
\newcommand{\Gr}{\text{Gr}}
\newcommand{\TrAd}{\text{TrAd}}

\newcommand{\Ad}{\text{Ad}}

\newcommand{\Gal}{\text{Gal}}

\newcommand{\Lie}{\text{Lie}}

\begin{document}

\title[characterization of complex Shimura varieties]%
{On the characterization\\ of complex Shimura varieties}
\author[Yakov Varshavsky]{Yakov Varshavsky}
\address{Department of Mathematics, University of Toronto,
100 St. George street, Toronto, Ontario M5S 3G3 CANADA}
\email{vyakov@math.toronto.edu}
\date{\today}
\maketitle
\begin{abstract}
In this paper we recall the construction and basic properties of complex
Shimura varieties and show that these properties actually
characterize them. This characterization immediately implies
the explicit form of Kazhdan's theorem on 
the conjugation of Shimura varieties. As a further corollary, 
we show that each Shimura variety corresponding to an
adjoint group has a canonical model over its reflex field.
We also indicate how this characterization implies the existence
of a $p$-adic uniformization of certain unitary Shimura varieties.
In the appendix we give a complete scheme-theoretic
proof of Weil's descent theorem.
\end{abstract}

\section{Introduction}

Shimura varieties are natural generalizations of classical modular curves.
As complex manifolds, Shimura varieties are just finite disjoint
unions of quotients of Hermitian symmetric domains by
arithmetic (congruence) subgroups. By the theorem of Baily-Borel they 
have canonical structures of complex quasi-projective varieties.

In some cases Shimura varieties have a modular interpretation as 
spaces parametrizing families of polarized abelian varieties with 
some additional structure,
or, more generally, of abelian motives. This interpretation have been
proven to be very useful for establishing various properties
of Shimura varieties, such as for showing that they
have the canonical structures (models) of  quasi-projective varieties defined over 
certain number fields, called the reflex fields, and for 
computing their cohomology.
 In particular, it enabled to show that zeta-functions 
of some simple Shimura varieties are products of automorphic 
$L$-functions, and that some  cases of Langlands' correspondence
can be realized in their cohomology.

Unfortunately not for all Shimura varieties there is known a 
modular interpretation. This difficulty was first overcomed by Kazhdan 
in \cite{Ka1, Ka2} (see also \cite{NR}), where he showed
that every conjugate of a Shimura variety is a Shimura
variety. This result combined with some  modular and group-theoretic methods allows
to prove the existence of the canonical models for all Shimura
varieties (see \cite{Brv, Mi1} and \cite{Mi3}).
The strategy of Kazhdan's proof was to
show that Shimura varieties are so rigid that the conjugate of a Shimura
variety can not be anything else than a Shimura variety. 
%
%

One of the main goals of this paper is to explain why Kazhdan's
theorem holds. For this we  show that complex Shimura varieties
together with a certain additional structure can be characterized by
their properties. Namely, we consider not just a Shimura variety but the triple
consisting of a Shimura variety $X$, a standard $G$-torsor $P$ 
with a flat connection $\C{H}$ over it, and an equivariant map $\rho$ from 
$P$ to a certain generalized Grassmannian $\Gr$. More
precisely, we consider an adelic family of the above triples, and
our Main Theorem (see Section 4) asserts that our family is the unique
family, satisfying certain group-theoretic properties.

The strategy of the proof is as follows: 
Let $(X;(P,\C{H});\rho:P\to\Gr)$ be the triple, and we want to show
that it is associated to a certain Shimura variety.
Let $M$ be the universal cover of some connected component of the
complex manifold $X\an$, and let $\Gm$ be its fundamental group. 
Then $X\an$ is isomorphic to $\Gm\bs M$, and our goal is to show that
certain properties of the triple imply that $M$ is a Hermitian
symmetric domain, and $\Gm\subset\Aut(M)$ is an arithmetic (congruence) subgroup. 

Using the flat connection $\C{H}$ on $P$, we conclude that the principal bundle
$P\an$
on $X\an$ is isomorphic to $\Gm\bs[M\times G\an]$. Then the composition of the
embedding $M\cong M\times\{1\}\hra P\an$ and $\rho$ gives us a
 $\Gm$-equivariant morphism $\rho_0:M\to \Gr\an$. Since our triple is
a part of adelic family, we conclude from our properties that
$\rho_0$ is a local isomorphism and that it is equivariant with respect to 
much larger group $\wt{\Gm}\supset\Gm$. Using the fact that $\Gr\an $ is a compact Hermitian
symmetric space, we conclude then from our properties that $M$ has to
be the Hermitian symmetric domain, dual to $\Gr\an$. 
Since $\Gm\bs M$ is a quasi-projective variety, we get that $\Gm$ is a lattice in
the semisimple Lie group $\Aut(M)$, and we conclude from our properties
that it is irreducible and has an infinite index in its
commensurator. Therefore by Margulis' theorem, $\Gm$ is arithmetic, as claimed. 

Since the properties, characterising the triple, are stable under the
conjugation, our Main Theorem immeditely implies the explicit
form of Kazhdan's theorem, which is a weak form of
Langlands' conjecture. Using the fact that Shimura varieties
corresponding to adjoint groups have no non-trivial equivariant
automorphisms, we then conclude that such Shimura varieties 
have unique equivariant structures over their reflex fields (which
coincide by uniqueness with the canonical models in the sense of Deligne).
For this we show that the descent datum, obtained from
Kazhdan's theorem is effective. This gives a completely uniform and group-theoretical 
construction of equivariant models over the reflex field in the case
of Shimura varieties, corresponding to adjoint groups.
However we do not know how to show directly that the constructed
models are canonical (that is, that the Galois group
acts on the set of the special points by Deligne's reciprocity law).
In particular, we do not know how to construct canonical models for 
arbitrary Shimura varieties.

We also indicate how our Main Theorem implies the existence 
the $p$-adic uniformization for some unitary Shimura varieties,
which was proven in \cite{Va1,Va2}.
For this we extend $p$-adically uniformized varieties to certain triples
and check that they satisfy all the properties, characterizing Shimura varieties.

The paper is organized as follows. In the second section we 
recall basic properties of Hermitian symmetric domains and of
locally symmetric varieties. In the third section we  recall how to
associate to a Shimura pair  an equivariant 
triple and show its basic properties. In the fourth section we formulate our 
Main Theorem and prove its consequences, while the fifth section 
is devoted to the proof of the Main Theorem. 
In the appendix we give a complete
scheme-theoretical proof of the classical Weil's descent theorem, 
which is crucially used in the construction of the canonical models of 
Shimura varieties.

 Finally we would like to remark that though for simplicity we work
only with Shimura varieties, corresponding to adjoint groups, 
most of the results and the proofs of the paper can be extended to the 
general case without almost no changes.

\centerline{\bf Notation and conventions}


1) For a Lie group $G$ let $G^0$ be its connected component of the identity.

2) For a totally disconnected topological group $E$ let 
$E^{disc}$ be the group $E$ with the discrete topology.




 
3) For an analytic space or a scheme $X$ (resp. an a point $x\in X$)
let $T(X)$ (resp. $T_x(X)$) be the tangent bundle on $X$ 
(resp. the tangent space at the point $x\in X$).



4) For a number field $F$ and a finite set $N$ of finite primes of
$F$ let $\B{A}_F^{\be}$ and $\B{A}_F^{\be,N}$ 
be the ring of finite adeles of $F$ and the ring of finite adeles of
$F$ without the components from $N$ respectively, and we omit
$F$ from the notation in the case $F=\B{Q}$. 
(Observe that this notation differs from that of \cite{Va1,Va2}).



5) For a scheme $X$ over a field $K$ and a $K$-algebra $A$ 
we will write $X_A$ or $X\otimes_K A$ instead of
$X\times_{\Spec K}\Spec A$.

6) For a complex analytic space $M$ let $\Aut(M)$ be the topological
group of holomorphic automorphisms of $M$
equipped with the compact-open topology.
 


 
7) For an algebraic variety $X$ over the field of complex numbers or 
over a $p$-adic field let $X\an$ the the 
corresponding complex or $p$-adic analytic space.

\centerline{\bf Acknowledgements}

First of all the author would like to thank R. Livn\'e, who stimulated
author's research on the 
$p$-adic uniformization of Shimura varieties \cite{Va1, Va2}, 
from which the ideas of the present paper were extracted.  

We also would like to thank O. Gabber, who noted to us, that the 
effectivity of the descent data of Shimura varieties is not obvious, 
to G. Faltings, who explained how to reformulate the continuity assumption of 
Weil's descent theorem into modern language, and to D. Blasius for his
valuable remarks.

This research was conceived while the author enjoyed the hospitality
and the financial support of the Max-Planck-Institute f\"ur Mathematik
in Bonn.
 
\section{Preliminaries} \label{S:loc}

We start from recalling some basic properties of (Hermitian) symmetric 
spaces of non-compact type, which we will call (Hermitian) symmetric 
domains, following \cite{He} and \cite[1.1-1.2]{De2}.

\begin{Emp} \label{E:Herm}
Let $\bo{G}$ be a real adjoint semi-simple group isotropic group, and
let $M^0$ be the set of all maximal compact subgroups of
$\bo{G}(\B{R})^0$. Then  the group $\bo{G}(\B{R})^0$ acts transitively
on $M^0$ (by conjugation), and the stabilizer in $\bo{G}(\B{R})^0$ 
of each point $K\in M^0$ is $K$. Moreover, $M^0$ has a
natural structure of a Riemannian symmetric domain such that the
connected component of the group of isomerties of $M^0$ is  
$\bo{G}(\B{R})^0$. 
Conversely, every Riemannian symmetric domain is obtained in this way
(see \cite[IV, V]{He}).
\end{Emp}

\begin{Emp} \label{E:conj}
It is well-known (see \cite[VIII, Thm. 6.1]{He}) that $M^0$ has a 
structure of a Hermitian symmetric domain if and only if 
each simple factor of each $K\in M^0$ either is a simple factor of  
$\bo{G}(\B{R})^0$ or has a non-discrete center. 
The last condition holds if and
only if there exists a continuous injective homomorphism 
$h$ from the unit circle $S^1:=\{z\in\B{C}\,|\,|z|<1\}$ to $\bo{G}(\B{R})$, whose centralizer in 
$\bo{G}(\B{R})^0$ is some $K\in M^0$ and which projects injectively to each
non-compact simple factor of $\bo{G}(\B{R})$.
We will call such a homomorphism {\em special}.

For each special homomorphism $h$, whose centralizer is $K$, 
there exists a unique $\bo{G}(\B{R})^0$-invariant
Hermitian structure on $M^0$ such that $Ad(h(z))$ induces the 
multiplication by $z$ on $T_K(M^0)$ for each $z\in S^1\subset\B{C}\m$. 
Moreover, this correspondence gives us a bijection between the set of
$\bo{G}(\B{R})^0$-conjugacy classes of special homomorphisms and 
the set of $\bo{G}(\B{R})^0$-invariant Hermitian structures on $M^0$.
In particular, the set of points of a Hermitian symmetric domain $M^0$ with 
$Aut(M^0)^0\cong\bo{G}(\B{R})^0$ is canonically
identified with elements of a certain $\bo{G}(\B{R})^0$-conjugacy class of 
special homomorphisms $h:S^1\hra\bo{G}(\B{R})$.
\end{Emp}

\begin{Emp} \label{E:inn}
For each special homomorphism $h:S^1\hra\bo{G}(\B{R})$, the involution 
$Inn(h(-1))$ defines an inner twisting between $\bo{G}$ and 
the unique compact real form of $\bo{G}_{\B{C}}$. Moreover, 
the isomorpism class of the twisting depends only on the conjugacy
class of $h$.
In particular, every $\bo{G}$ for which a special homomorphism 
$h$ exists contains a maximal elliptic torus
that is, $\bo{G}(\B{R})$ contains a maximal compact torus of 
$\bo{G}(\B{C})$.
\end{Emp}

\begin{Emp} \label{E:mu}
Each special homomorphism $h$ from \re{conj} extends uniquely to an 
injective cocharacter $\mu=\mu_h:\B{G}_m\to\bo{G}_{\B{C}}$. 
By the properties of the Cartan decomposition, $\mu$ is {\em special}
in the following sense: the representation $Ad\circ\mu$  of $\B{G}_m$ on  
$Lie\bo{G}_{\B{C}}$ has no other weights than  $-1,0$ and $-1$
(compare \cite[Cor. 1.1.17 and Prop. 1.2.2]{De2}). 
Conversely, each special cocharacter $\mu$ is obtained in this way: 
choose a real compact form $\bo{G}^{comp}$ of $\bo{G}_{\B{C}}$ such that 
$\mu(S^1)\subset\bo{G}^{comp}(\B{R})$, let $\bo{G'}$ be the inner form
of  $\bo{G}^{comp}$, obtained by the involution $Inn(\mu(-1))$, and
let $h':S^1\to\bo{G'}(\B{R})$ be the restriction of $\mu$ to
$S^1$. Then $h'$ is a special homomorphism (compare \cite{De2}).

Thus the correspondence $(\bo{G'},h')\mapsto\mu$ gives a
bijection between the set of $\bo{G}(\B{C})$-conjugacy classes of 
pairs consisting of a real form $\bo{G'}$ of $\bo{G}_{\B{C}}$ and a 
special homomorphism $h':S^1\hra\bo{G'}(\B{R})$ and
the set of $\bo{G}_{\B{C}}$-conjugacy classes of special cocharacters of
$\bo{G}_{\B{C}}$.

Denote the map $\B{C}\m\to S^1\;(z\mapsto z/\overline{z})$ by $\pi$.
Then for each special homomorphism $h$, the composition
$\pi\circ h:\B{C}\m\to\bo{G}(\B{R})$ satisfies Deligne's axioms
(\cite[1.5]{De1}). Conversely, every homomorphism 
$\B{C}\m\to\bo{G}(\B{R})$, satisfying Deligne's axioms, factors as a 
composition of $\pi$ and a certain special homomorphism $S^1\to\bo{G}(\B{R})$.
\end{Emp}

\begin{Emp} \label{E:dual}
For each special cocharacter $\mu$ (see \re{mu}) and each $i=-1,0,1$, 
let $\Lie_i$ be the
weight $i$ subspace of $\Lie\,\bo{G}_{\B{C}}$ with respect to the action of 
$\B{G}_m$. Then the sum $\Lie_0\oplus\Lie_{-1}$ is 
a Lie algebra of a certain parabolic subgroup
$\bo{P}_{\mu}\subset\bo{G}_{\B{C}}$. Let $\C{D}_\mu$ be the subset of
the Dynkin diagram $\C{D}$ of $\bo{G}_{\B{C}}$, corresponding to
$\bo{P}_{\mu}$.
Then the parabolic subgroup $\bo{P}_{\mu}$ and a subset $\C{D}_{\mu}$
are {\em special} in the following sense (compare \cite[1.2.5]{De2}):

i) $\C{D}_{\mu}$ contains at most one vertex of each connected component
of $\C{D}$;

ii) each vertex of $\C{D}_{\mu}$ appears with multiplicity one in the 
decomposition of the highest root of $\bo{G}_{\B{C}}$.

Moreover, the correspondence 
$\mu\mapsto\bo{P}_{\mu}\mapsto\C{D}_{\mu}$ induces a bijection
between the set of conjugacy classes of special cocharacters 
of $\bo{G}_{\B{C}}$ , the set of conjugacy classes
of special parabolics of $\bo{G}_{\B{C}}$, and the set of special
subsets of the Dynkin diagram of $\bo{G}_{\B{C}}$.
\end{Emp}

\begin{Emp} \label{E:grass}
For each $h\in M^0$, as in \re{conj}, let
$\check{M}_{\mu_h}$ be the set of all
parabolic subgroups of $\bo{G}_{\B{C}}$, which are conjugate to 
$\bo{P}_{\mu_h}$. Then  
$\check{M}_{\mu_h}$ does not depend on $h$, and
$\check{M}:=\check{M}_{\mu_h}\cong\bo{G}(\B{C})/\bo{P}_{\mu_h}(\B{C})$ 
has a natural structure of a complex projective variety.
We will call $\check{M}$ the generalized Grassmannian, corresponding to
$(\bo{G},M)$.

The map $h\mapsto\bo{P}_{\mu_h}$ defines a holomorphic 
$\bo{G}(\B{R})^0$-equivariant map $\rho_0:M^0\to\check{M}\an$.
We claim that $\rho_0$ is the unique $\bo{G}(\B{R})^0$-equivariant
holomorphic map from $M^0$ to $\check{M}\an$. 

Indeed, let  $\rho'$ be any such a map, let $h$ be an element of $M^0$, 
and let $K\subset\bo{G}(\B{R})^0$ be the centralizer of $h$.
Then $\rho'(h)$ is a parabolic subgroup of $\bo{G}_{\B{C}}$, and the 
differential
of $\rho'$ at $h\in M^0$ induces a $\Lie\bo{G}$-equivariant $\B{C}$-linear
map 
$d\rho'_h:T_h(M^0)=Lie\bo{G}/Lie\,K\to T_{\rho'(h)}(\check{M}\an)=
\Lie\bo{G}_{\B{C}}/Lie\,\rho'(h)$.
By the properties of the Cartan decomposition, $d\rho'_h$ has to be an
 isomorphism.
Hence $\Lie\,\rho'(h)$ is  the kernel of the projection
 $\Lie\bo{G}_{\B{C}}\to T_h(M^0)$, 
obtained
as a unique $\B{C}$-linear extension of the natural projection
$\Lie\bo{G}\to T_h(M^0)$. This shows that $\rho'(h)$ equals 
$\bo{P}_{\mu_h}=\rho_0(h)$, showing the uniqueness property of $\rho_0$.

Moreover, the above argument also shows that $\rho_0$ is an open embedding
(which is called the Harish-Chandra-Borel embedding) and that 
$\check{M}\an$ is the compact dual of $M^0$ in the sense of 
\cite[V, $\S$2]{He}.
\end{Emp}


\begin{Emp} \label{E:hodge}
For every $h\in M^0$ and every real finite-dimensional representation 
$\tau:\bo{G}\to\Aut(V)$ the composition $\tau\circ h\circ\pi$, where 
$\pi$ is as in \re{mu}, defines
a Hodge structure on $V$ (see \cite[1.1.1]{De2}). 
Thus for every fixed $\tau$ we get a $\bo{G}(\B{R})^0$-equivariant
continuous family of Hodge structures on $V$, parametrized by $M^0$.
Moreover, \cite[Prop. 1.1.14]{De2} shows that this family is a
variation of Hodge structures over $M^0$.
\end{Emp}

\begin{Emp} \label{E:triple}
Let $\Gm\subset{\bo{G}}(\B{R})^0$ be a torsion-free lattice.  Then
we can form a triple consisting of:

i) a smooth complex manifold $\wt{X}_{\Gm}:=\Gm\bs M^0$;

ii) a principal $(\bo{G}_{\B{C}})\an$-bundle 
$\wt{P}_{\Gm}:=\Gm\bs[M^0\times(\bo{G}_{\B{C}})\an]$ on $\wt{X}_{\Gm}$ 
(where $\Gm$ acts on 
$M\times(\bo{G}_{\B{C}})\an$ by the formula
$\gm(m,g)=(\gm(m),g\gm^{-1})$),
equipped  with a natural flat connection
$\wt{\C{H}}\subset T(\wt{P}_{\Gm})$, consisting of vectors, tangent to
$M^0$;

iii) a $\bo{G}(\B{C})$-equivariant map 
$\wt{\rho}:\wt{P}_{\Gm}\to\check{M}\an$ defined by 
$\wt{\rho}(m,g)=g(\rho_0(m))$ for all $m\in M^0$ and $g\in \bo{G}_{\B{C}}$.

\end{Emp}

Recall the following important result.

\begin{Thm}  \label{T:BB}
Suppose that  $\Gm$ is arithmetic. Then:

a) The complex manifold $\wt{X}_{\Gm}$ has a unique 
structure $X_{\Gm}$ of a complex quasi-projective variety. Moreover,
if $Y$ is a complex algebraic variety, then every analytic map from
$Y\an\to\wt{X}_{\Gm}$ comes from an algebraic map from $Y$ to 
$X_{\Gm}$.

b) The pair $(\wt{P}_{\Gm},\wt{\C{H}})$ has a unique algebraic
structure of
a $\bo{G}_{\B{C}}$-torsor $P_{\Gm}$ on $X_{\Gm}$ equipped with a flat connection $\C{H}$.

c) The map $\wt{\rho}$ comes from an algebraic map $\rho$.
\end{Thm}
\begin{pf}
a) Since every hermitian symmetric domain is isomorphic to a bounded
 symmetric domain (see \cite[VIII, Thm. 7.1]{He}), the existence in
 the first statement of the theorem follows from (\cite{BB}),
 and the second statement follows from (\cite{Bo2}).

b) Deligne proved (see \cite[IV]{Bo3}) that for every smooth complex
algebraic variety $X$ the functor $(V,\nabla)\mapsto(V\an,\nabla\an)$
defines an equivalence between the category of vector bundles with
regular flat connections on $X$ and that of vector bundles with
flat connections on $X\an$. By Tannakian arguments,
 Deligne's theorem also holds for principal bundles
(compare \cite[Thm. 3.2]{DM} or \cite[Prop.1.9.13]{Va1}). 

c) Let $k$ be the complex dimension of $M^0$. Then the map
$i:\bo{P}\mapsto\Lie\bo{P}\subset\Lie\bo{G}_{\B{C}}$ defines a closed 
embedding of $\check{M}$ into Grassmannian
$\Gr_k:=\Gr_k(\Lie\bo{G}_{\B{C}})$, parametrizing the set of 
all linear subspaces of $\Lie\bo{G}_{\B{C}}$ of codimension $k$.
Our statement will follow if we show that the composition map 
$i\circ\wt{\rho}:(P_{\Gm})\an\to\Gr\an$  
is algebraic. By the universal property of Grassmannians, we have to
check that $V:=\{[p,v]\in P_{\Gm}\times\Lie\bo{G}_{\B{C}}]\,|\,
v\in\Lie\,\wt{\rho}(p)\}$ is an algebraic subbundle of 
the trivial bundle $P_{\Gm}\times\Lie\bo{G}_{\B{C}}$  on $P_{\Gm}$.

By the descent theory (see \cite[I, 2.21-2.22]{Mi2} or 
\cite[Lem. 1.9.5]{Va1}), the quotient $V_{\Gm}$ of  
$P_{\Gm}\times\Lie\bo{G}_{\B{C}}$ by the diagonal action of
$\bo{G}_{\B{C}}$ is a vector bundle on $X_{\Gm}$.
Moreover, $(V_{\Gm})\an\cong\Gm\bs[M\times(\Lie\bo{G}_{\B{C}})\an]$, 
where the quotient is taken with respect to the
diagonal action of $\Gm$. By \re{hodge},
$(V_{\Gm})\an$ has a natural structure of a variation of Hodge
structures, hence by the  theorem of Griffiths
(see \cite[Thm. 4.13]{Sc}) every holomorphic subbundle, which
lies in the Hodge filtration of $(V_{\Gm})\an$, 
is an algebraic subbundle of $V_{\Gm}$. Since
$V\subset P_{\Gm}\times\Lie\bo{G}_{\B{C}}$ is the inverse image
of the  zero term of the 
Hodge filtration of $(V_{\Gm})\an$, it is therefore an algebraic subbundle, 
as was claimed.
\end{pf}

\section{Basic construction}

\begin{Emp} \label{E:data}
We start from a pair consisting of a simple
adjoint group $\bo{G}$ over $\B{Q}$, 
and a $\bo{G}(\B{R})$-conjugacy class $M$ of special homomorpisms
$h:S^1\hra\bo{G}(\B{R})$, as in \re{conj}, which we will call {\em
Shimura pair}.
 By the remark at the end of \re{mu}, our notion of Shimura pair
essentially coinsides with that of \cite[1.5]{De1}.

By \re{conj}, $M$ has a natural structure of a complex
manifold such that each connected component of $M$ 
is a Hermitian symmetric domain. Moreover, to the pair 
$(\bo{G},M)$ we associate (as in  \re{dual}) a generalized 
Grassmannian $\check{M}$, and we 
have a unique $\bo{G}(\B{R})$-equivariant holomorphic map 
$\rho_0:M\to\check{M}\an$, which is an open embedding.
Notice also that the uniqueness of $\rho_0$ implies that $M$ has no 
non-trivial $\bo{G}(\B{R})$-equivariant holomorphic automorphisms.
\end{Emp}

Next we recall how to associate to the pair $(\bo{G},M)$ a certain 
$\bo{G}(\B{A}^{\be})$-equivariant triple, which is basic for this paper.


\begin{Emp} \label{E:Shim}
Consider the complex analytic space 
$\wt{X}=[M\times\bo{G}(\B{A}^{\be})^{disc}]/\bo{G}(\B{Q})$,
where the group $\bo{G}(\B{Q})$ acts on the product 
$M\times\bo{G}(\B{A}^{\be})$ by the action $(m,g)\gm=(\gm^{-1}(m),g\gm)$.
Then the group $\bo{G}(\B{A}^{\be})$ acts analytically on $\wt{X}$ by the left
multiplication.

Set $\bo{G}(\B{Q})_+:=\bo{G}(\B{Q})\cap\bo{G}(\B{R})^0$, and put 
$\Gm_T:=\bo{G}(\B{Q})_+\cap T$ for every compact and open subgroup 
$T\subset\bo{G}(\B{A}^{\be})$. Then for every compact and open subgroup 
$S\subset\bo{G}(\B{A}^{\be})$ the quotient $S\bs\wt{X}$ is a finite
disjoint union of the quotients ${\Gm_{aSa^{-1}}}\bs M^0$, where $M^0$ is 
a connected component of $M$ and
$a$ runs over the set of representatives of the double classes
$S\bs\bo{G}(\B{A}^{\be})/\bo{G}(\B{Q})_+$. 
Since each  $\Gm_T$ is an arithmetic congruence subgroup, \rt{BB}, a) 
implies that  $\wt{X}_S$ has a unique structure $X_S$ of
a complex smooth quasi-projective scheme. Moreover, for all $T\subset S$ the natural 
analytic projection $\wt{X}_T\to\wt{X}_S$ comes from a finite 
 algebraic map $X_T\to X_S$, which is etale
 for all sufficiently small $S$.
Hence the $X_S$'s form a projective system, for which all the
transition maps are affine. Therefore there exists an inverse limit 
$X=X(\bo{G},M)$ of the $X_S$'s in the 
category of schemes, which we will call the complex {\em Shimura variety}
corresponding to Shimura pair $ (\bo{G},M)$.

Notice next that $\wt{X}$ is the inverse limit of the $X_S\an$'s 
in the category of the analytic spaces and that the
transitions maps  $X_T\an\to X_S\an$ are topological coverings
for all sufficiently small $T\subset S$.
Hence $X\an:=\wt{X}$ satisfies the universal property of an analytic 
space, corresponding to $X$ (see \cite[Exp. XII, Thm. 1.1]{SGA1}). 
By the second statement of \rt{BB}, a), the action
of $\bo{G}(\B{A}^{\be})$ on $\wt{X}$ defines  an action of 
$\bo{G}(\B{A}^{\be})$ on $X$. Then the quotient $S\bs X$
is canonically isomorphic to $X_S$ for every compact and open 
$S\subset\bo{G}(\B{A}^{\be})$. 
\end{Emp}

\begin{Emp}
Similarly we can construct (using \rt{BB}, b) and c)) a 
$\bo{G}(\B{A}^{\be})$-equivariant $\bo{G}_{\B{C}}$-torsor $\pi:P\to X$
equipped with a flat connection $\C{H}$ and a  
$\bo{G}(\B{A}^{\be})\times\bo{G}(\B{C})$-equivariant 
map $\rho:P\to\check{M}$ such that 
$P\an$ is isomorphic to 
$[M\times\bo{G}_{\B{C}}\an\times\bo{G}(\B{A}^{\be})^{disc}]/\bo{G}(\B{Q})$,
equipped with a natural connection
and $\rho\an[m,g_{\B{C}},g_f]=g_{\B{C}}(\rho_0(m))$ for all 
$m\in M$, $g_{\B{C}}\in\bo{G}_{\B{C}}$ and 
$g_f\in\bo{G}(\B{A}^{\be})$ (compare \cite[4.2 and Prop. 4.3.2]{Va1}).
Since $\rho_0$ is an open embedding, the differential $d\rho$ induces
an isomorphism between $\C{H}_p$ and $T_{\rho(p)}(\check{M})$ 
for each $p\in P$.
\end{Emp}

\begin{Emp} \label{E:shimtr}
{\bf Summary.} In the previous subsections we have associated to 
Shimura pair $(\bo{G},M)$ a 
$\bo{G}(\B{A}^{\be})$-equivariant triple $T=T(\bo{G},M)$ consisting of:

i) A scheme $X$ over $\B{C}$ such that for each compact and open subgroup $S\subset\bo{G}(\B{A}^{\be})$
the quotient $X_S:=S\bs X$ exists and is a complex smooth quasi-projective scheme.  Moreover, $X$ is the inverse limit of the $X_S$'s, and
the transition maps $X_{S_1}\to X_{S_2}$ are finite and etale for all 
sufficiently small $S_1\subset S_2$. 

ii) a $\bo{G}_{\B{C}}$-torsor $P$ on $X$ equipped with a flat
 connection $\C{H}$.

iii) a $\bo{G}_{\B{C}}$-equivariant map $\rho:P\to\check{M}$ such that 
its differential $d\rho$ induces an isomorphism between $\C{H}_p$ and 
$T_{\rho(p)}(\check{M})$ for each $p\in P$.
\end{Emp}

We will call $T(\bo{G},M)$ the triple corresponding to the Shimura
variety of $(\bo{G},M)$.

\begin{Rem} \label{R:isom}
a) Every morphism $\varphi:(\bo{G_1},M_1)\to(\bo{G_2},M_2)$ of two
Shimura pairs (that is, $\varphi$ is a group homomorphism
$\bo{G_1}\to\bo{G_2}$, which maps $M_1$ to $M_2$) defines (by \rt{BB} and its proof)
a $\varphi$-equivariant morphism $T(\bo{G_1},M_1)\to
T(\bo{G_2},M_2)$ of the corresponding triples 
(compare \cite[Lem. 2.2.6 and Rem. 4.4.4]{Va1}). 
 In particular, the isomorphism class of
the triple $T(\bo{G},M)$ is determined by the isomorphism class of
the pair $(\bo{G},M)$.

b) Conversely, the isomorphism class of the triple $T(\bo{G},M)$, 
and even the isomorphism class of the Shimura variety $X(\bo{G},M)$ 
determines the isomorphism class of the pair $(\bo{G},M)$.
\begin{pf}
Suppose that we have a continuous isomorpism
$\varphi:\bo{G_1}(\B{A}^{\be})\isom\bo{G_2}(\B{A}^{\be})$ and a 
$\varphi$-equivariant isomorphism 
$\psi:X_1=X(\bo{G_1},M_1)\isom X_2=X(\bo{G_2},M_2)$. Then $\psi$
induces the analytic isomorphism 
$\psi\an:X_1\an=[M_1\times\bo{G_1}(\B{A}^{\be})^{disc}]/\bo{G_1}(\B{Q})\isom
X_2\an=[M_2\times\bo{G_2}(\B{A}^{\be})^{disc}]/\bo{G_2}(\B{Q})$.
Hence there exists an element $g\in \bo{G_2}(\B{A}^{\be})$ and connected 
components $M_1^0$ and $M_2^0$ of $M_1$ and $M_2$ respectively such that
$\psi\an$ maps the connected component $M_1^0\times\{1\}$ of
$X_1\an$ into $M_2^0\times\{g\}$. After replacing
$\varphi$ by $Inn(g)\circ\varphi$ and $\psi$ by $g\circ\psi$, 
we may assume that $g=1$. Then $\varphi$ induces
an isomorphism  $\bo{G_1}(\B{Q})_+\isom\bo{G_2}(\B{Q})_+$ 
between the stabilizers $\bo{G_1}(\B{Q})_+$ and 
$\bo{G_2}(\B{Q})_+$ of $M_1^0\times\{1\}$ and $M_2^0\times\{1\}$
respectively. Since $\varphi$ is continuous, and since each 
$\bo{G_i}(\B{Q})_+$ is Zariski dense in $\bo{G_i}$, we get that
 $\varphi$ is induced by certain algebraic isomorphism $\varphi_0$ between 
$\bo{G_1}$ and $\bo{G_2}$. Since each $\bo{G_i}(\B{Q})$ is dense in 
$\bo{G_i}(\B{R})$, the map  $\psi\an$ induces a 
$\varphi_0$-equivariant biholomorphic map
$\psi_0:M_1\cong M_1\times\{1\}\isom M_2\times\{1\}\cong M_2$. The
last observation in \re{data} now implies that $\psi_0$ has to
coincide with the map, induced by $\varphi_0$. This shows the assertion.
\end{pf}
\end{Rem} 

For each $p\in P(\B{C})$ let 
$\Gm_p\subset\bo{G}(\B{C})\times\bo{G}(\B{A}^{\be})$ be the stabilizer of $p$.

\begin{Prop} \label{P:prop}
The triple $T(\bo{G},M)$ satisfies the following properties:

a) The group $\bo{G}(\B{A}^{\be})$ acts transitively on the set of the 
connected components of $X$.

b) The $\bo{G}(\B{A}^{\be})$-structure group of 
$(P,\C{H})$ can not be $\bo{G}(\B{A}^{\be})$-equivariantly reduced to
a proper subgroup of $\bo{G}_{\B{C}}$.

c) If $S$ is sufficiently small, $X_S$ does not contain a nontrivial 
homogeneous space for the action of $\bo{G}$.

d) There exists a point $p\in P$ such that 
the Zariski closure of the projection of $\pr_{\B{C}}(\Gm_p)\subset\bo{G}(\B{C})$ to $\Aut(\check{M})$ contains a maximal torus.

e) Each  $\gm\in\pr_{\B{C}}(\Gm_p)$ is semisimple, 
and all the eigenvalues of $\Ad\,\gm$ have absolute value one.

f) Each $\Gm_p$ is ``diagonally embedded'' into
$\bo{G}(\B{C})\times\bo{G}(\B{A}^{\be})$ ``up to a conjugation'' 
in the following sense:
for each number field $L$ and each conjugation-invariant regular function
$f$ on $\bo{G_L}$, we have 
$f(\Gm_p)\subset L$, embedded diagonally into 
$L\otimes_{\B{Q}}(\B{C}\times\B{A}^{\be})$.

g) Each group $\Gm_p$ is isomorphic to the group of $\B{Q}$-rational
points of a certain algebraic group $\bo{G_p}$ over $\B{Q}$.
\end{Prop}
\begin{pf}

a) follows, since the group $\bo{G}(\B{A}^{\be})$ acts transitively on the 
set of the analytic connected components of $X\an$. 

b) follows from (and is actually equivalent to) the fact that 
$\bo{G}(\B{Q})_+$ is Zariski dense in $\bo{G}_{\B{C}}$ (compare
\cite[Prop. 4.5.1]{Va1}). Indeed, if $\bo{G'}$ is a proper algebraic subgroup
of $\bo{G}_{\B{C}}$ and if $P'\subset P$ is a  
$\bo{G}(\B{A}^{\be})$-equivariant $\bo{G'}$-torsor over
$X$ with a flat connection such that the natural map 
$\bo{G}_{\B{C}}\times_{\bo{G'}}P'\to P$ is an isomorphism of torsors,
preserving the connections, then there exists an element
$g\in\bo{G}(\B{C})$ such that 
$(P')\an\cong [M^0\times(\bo{G'})\an g\times\bo{G}(\B{A}^{\be})^{disc}]/
\bo{G}(\B{Q})_+$ for a certain connected component $M^0$ of $M$.
Then $g\bo{G}(\B{Q})_+ g^{-1}\subset \bo{G'}$, hence the Zariski
density of $\bo{G}(\B{Q})_+$ in $\bo{G}_{\B{C}}$ implies that 
$\bo{G'}=\bo{G}_{\B{C}}$, as was claimed.

c) If $S$ is sufficiently small, the universal cover of each 
connected component of $X_S$ is a Hermitian symmetric domain. 
Therefore $X_S\an$ is a Kobayashi hyperbolic space, hence any action of a 
linear algebraic groups on $X_S$ has to be trivial
(see, for example, \cite[I, $\S$2]{La}).

d) By \re{conj}, the group $\bo{G}_{\B{R}}$ contains a maximal elliptic torus 
$\bo{T}$.
Therefore the set of regular elliptic elements
of $\bo{G}(\B{R})^0$ (that is, of those elements of $\bo{G}(\B{R})^0$ whose
centralizer is a maximal elliptic torus) is open, because it is equal to the
union of all the conjugates of the regular elements in
$\bo{T}(\B{R})$. Since $\bo{G}(\B{Q})_+$
is dense in $\bo{G}(\B{R})^0$, it contains some regular elliptic element 
$\dt$. Then the centralizer $\bo{T_{\dt}}\subset\bo{G}$ is a maximal
$\B{Q}$-rational torus such that $\bo{T_{\dt}}(\B{R})$ is compact.
In particular, $\bo{T_{\dt}}(\B{R})\subset\bo{G}(\B{R})^0$ stabilizes 
some point $m\in M$. Then the stabilizer of the point 
$p:=[m,1,1]\in P(\B{C})$ contains 
$\bo{T_{\dt}}(\B{Q})\subset\bo{G}(\B{Q})_+\subset\bo{G}(\B{C})\times\bo{G}(\B{A}^{\be})$, 
implying the assertion (compare also \cite[Thm. 5.1]{De1}).

e)-g) For each $g\in\bo{G}(\B{C})\times\bo{G}(\B{A}^{\be})$ and each
$p\in P$ we have $\Gm_{g(p)}=g\Gm_p g^{-1}$. Since the properties e)-g) 
are stable under the conjugation, we may therefore assume that 
$p$ is of the form $[m,1,1]\in M\times\{1\}\times\{1\}\subset
P(\B{C})$, hence $\Gm_p=Stab_{\bo{G}(\B{Q})}(m)$ for some $m\in M$.
For such a point

e) holds, because the stabilizer of every point in $M$ is compact;

f) holds, because every $\Gm_p$ is contained in $\bo{G}(\B{Q})$, 
embedded diagonally;

g) holds (with $\bo{G_p}$ equal to the Zariski closure of $\Gm_p$ in $\bo{G}$), 
because the stabilizer of every $m\in M\subset\check{M}(\B{C})$ in $\bo{G}$  
is an algebraic subgroup.
\end{pf}

\begin{Not}
Let $F$ and $\bo{H}$ be the number field and the 
absolutely simple adjoint group over it such that $\bo{G}$ 
is the restriction of scalars of $\bo{H}$ from $F$ to $\B{Q}$.
\end{Not}

\begin{Rem} \label{R:propf}
Property f) of \rp{prop} is equivalent to the property 

f') For each $p\in P$, each number field $L$ and each $L$-rational representation 
$\mu$ of $\bo{G_L}$, we have 
$Tr\mu(\Gm_p)\subset L\subset L\otimes_{\B{Q}}(\B{C}\times\B{A}^{\be})$, 
embedded diagonally.

In particular, f) implies the following properties:

 f${}_1$) For each $p\in P$ we have 
$\TrAd_{\bo{H}}(\Gm_p)\subset F\subset F\otimes_{\B{Q}}
(\B{C}\times\B{A}^{\be})$, embedded diagonally. 
(Here $\Gm_p$ is considered as a subgroup of 
$\bo{H}(F\otimes_{\B{Q}}\B{C})\times\bo{H}(\B{A}^{\be}_F)$, and $Ad_{\bo{H}}$
is the adjoint representation of $\bo{H}$.)  

 f${}_2$) For each $p\in P$ we have 
$\TrAd_{\bo{G}}(\Gm_p)\subset \B{Q}\subset \B{C}\times\B{A}^{\be}$, 
embedded diagonally.

If we assume the property e), then  f${}_1$) in its turn implies the
following properties:

 f${}_3$) Each $\Gm_p$ projects injectively to every 
simple factor of $\bo{G}(\B{C})\times\bo{G}(\B{A}^{\be})$.

 f${}_4$) The group $\bo{G}(\B{A}^{\be})$ acts faithfully on $X$.
\end{Rem}
\begin{pf}
The equivalence of f) and f') follows from the fact that the algebra 
of conjugation-invariant regular functions of a reductive group
 is generated as a vector space by the traces of irreducible 
representations. 

Taking $\mu$ be the adjoint representations of
$\bo{H}$ (considered as a representation of $\bo{G}$ defined over 
$F$) and the adjoint representation of $\bo{G}$, we get  
f${}_1$) and f${}_2$) respectively.

By e), the group $\pr_{\B{C}}(\Gm_p)$ consists of semi-simple elements only,
therefore f${}_1$) implies  f${}_3$). 

Finally, the kernel of the action of $\bo{G}(\B{A}^{\be})$ on $X$ is a 
normal subgroup. Hence it has to be a product of normal subgroups 
of simple factors of $\bo{G}(\B{A}^{\be})$, thus it has to be trivial 
by f${}_3$).
\end{pf}

\begin{Rem} \label{R:propg} 
 Properties e) and f${}_2$) of \rp{prop} 
imply that each $\Gm_p$ can be naturally embedded as a Zariski dense 
subgroup of $\bo{G_p}(\B{Q})$ for a certain reductive (not necessary connected) group 
$\bo{G_p}$ over $\B{Q}$. Thus property g) just says 
that this embedding has to be surjective.
\end{Rem}
\begin{pf}
Let $\bo{P_p}$ be the Zariski closure of $\pr_{\B{C}}(\Gm_p)$ in 
$\bo{G}_{\B{C}}$, let $\bo{M_p}$ be a Levi subgroup of  
$\bo{P_p}$, and let 
$\pr':\Gm_p\to\bo{M_p}(\B{C})\subset\bo{G}(\B{C})$ be the 
composition of $\pr_{\B{C}}$ with the natural projection 
$\bo{P_p}\to\bo{M_p}$. Properties e) and f${}_2$) imply,
as in \rr{propf}, that the map $\pr'$ is injective. 
Moreover, property  f${}_2$) implies that 
$\TrAd_{\bo{G}}\pr'(\Gm_p)=\TrAd_{\bo{G}}\pr_{\B{C}}(\Gm_p)\subset\B{Q}$.

Now we argue as in the proof of \cite[Ch. VIII, Prop. 3.20]{Ma}:
Let $M$ and $M_0$ be the $\B{C}$-span of 
$Ad_{\bo{G}}(\bo{M_p}(\B{C}))\subset End(\Lie\bo{G}_{\B{C}})$ and the
$\B{Q}$-span of 
$Ad_{\bo{G}}(\pr'(\Gm_p))\subset End(\Lie\bo{G}_{\B{C}})$
respectively. Since $\bo{M_p}$ is reductive, 
$\Lie\bo{G}_{\B{C}}$ is a semisimple representation of $\bo{M_p}$.
Then we deduce from Jacobson's density theorem 
(see \cite[XVII, $\S$3]{La1}) that $M_0$ is a $\B{Q}$-structure of $M$.

Let the group $Aut(\Lie\bo{G}_{\B{C}})$ acts on 
$End(\Lie\bo{G}_{\B{C}})$ by the left multiplication. Then 
$\Gm_p\cong\Ad_{\bo{G}}(\pr'(\Gm_p))$ maps $M_0$ into itself, and 
we define the required group $\bo{G_p}$ be the Zariski closure of the 
image of $\Gm_p$ in $\bo{GL}(M_0)$.
\end{pf}

\begin{Lem} \label{L:aut}
We have the following properties:

1) The group of automorphisms of $X_S$ is finite for each 
sufficiently small $S$.

2) There is no non-trivial  $\bo{G}(\B{A}^{\be})$-equivariant automorphism
 of $X$.

3) There is no non-trivial 
$\bo{G}(\B{A}^{\be})\times\bo{G}_{\B{C}}$-equivariant automorphisms of $P$ 
over $X$, preserving $\C{H}$.

4) $\rho$ is the unique  
$\bo{G}_{\B{C}}\times\bo{G}(\B{A}^{\be})$-equivariant map from 
$P$ to $\check{M}$.
\end{Lem}
\begin{pf}
1) Since $X_S$ has only finitely many connected components,
it is enough to show the finiteness of the automorphism group
of each connected component $X_S^0$ of $X_S$. 
But if $S$ is sufficiently small, then each 
$(X_S^0)\an$ is isomorphic to $\Gm\bs M^0$ for some torsion-free congruence subgroup
$\Gm\subset\bo{G}(\B{Q})_+\subset\bo{G}(\B{R})^0$. 
Hence $Aut(X_S^0)\cong Norm(\Gm)/\Gm$,
where by $Norm(\Gm)$ we denote the normalizer of $\Gm$ in
$\bo{G}(\B{R})^0$. 

Since  $\Gm$ is Zariski dense in $\bo{G}$ (see, for
example, \cite[Ch. 4, Thm. 4.10]{PR}), it
has a trivial centralizer in $\bo{G}(\B{R})^0$.
Since $\Gm$ is discrete in $\bo{G}(\B{R})^0$, we therefore see
that $Norm(\Gm)$ has to be discrete as well.
But $\Gm$ has a finite covolume in $\bo{G}(\B{R})^0$ (compare 
Step 2, I) of \rs{proof}), hence it
has to be of finite index in  $Norm(\Gm)$, as was claimed.

2) The proof of \rr{isom}, b) shows that every 
$\bo{G}(\B{A}^{\be})$-equivariant automorphism of $X$
has to be induced by an element of $\bo{G}(\B{A}^{\be})$. Since $\bo{G}$
is adjoint, this element has therefore to be the identity.

3) Any $\bo{G}(\B{C})\times\bo{G}(\B{A}^{\be})$-equivariant automorphism 
of $(P,\C{H})$ induces a $\bo{G}(\B{Q})$-equivariant
automorphism $\psi$ of the trivial torsor $M\times\bo{G}_{\B{C}}$ over 
$M$, preserving the trivial connection. Hence there exists a holomorphic map
$\phi:M\to(\bo{G}_{\B{C}})\an$ such that $\psi(m,g)=(m,g\phi(m))$
for all $m\in M$ and $g\in(\bo{G}_{\B{C}})\an$.
Since $\psi$ preserves the trivial connection, $\phi$ is a
constant map. But $\psi$ is $\bo{G}(\B{Q})$-equivariant, hence $\phi$ has then to
be trivial, as was claimed.

4) Let $\rho'$ be any map, as in the statement. Then the restriction of 
$(\rho')\an$ to $M\cong M\times\{1\}\times\{1\}\subset P\an$ is a 
$\bo{G}(\B{Q})$-equivariant holomorphic map 
$\rho'_0:M\to\check{M}\an$. Since $\bo{G}(\B{Q})$ is dense in 
$\bo{G}(\B{R})$, the map $\rho'_0$ is $\bo{G}(\B{R})$-equivariant,
so it has to coincide with $\rho_0$ (see \re{data}). Since  
$\rho'$ is $\bo{G}_{\B{C}}\times\bo{G}(\B{A}^{\be})$-equivariant,
it has to coincide therefore with $\rho$.  
\end{pf}

\section{Main Theorem and its applications} \label{S:Main}

\begin{Main} \label{T:Main}
Suppose that we are given a Shimura pair $(\bo{G},M)$ and a triple $T$ 
consisting of data i)-iii) of \re{shimtr} and  satisfying the 
properties a)-g) of \rp{prop}. 
Then $T$ is isomorphic to the triple $T(\bo{G},M)$, corresponding to the 
Shimura variety of $(\bo{G},M)$. 
\end{Main}

\begin{Rem} \label{R:prop}
a) Not all of the properties a)-g) of \rp{prop} are equally important, 
and if we omit or weaken some of them, then $X$ (resp. $T$) 
will only slightly differ from the Shimura variety (resp. from the
triple), corresponding to $(\bo{G},M)$.
Namely small modifications of the proof of the Main Theorem
(see \rs{proof}) show that:

i) We need property a) only to assure that $X$ is a Shimura
variety and not a finite disjoint unions of Shimura varieties.

ii) We need property g) only to assure that  $X$ is a Shimura variety and 
not an equivariant finite etale cover of it. In particular, g)
automatically holds if the group $\bo{G}(\B{Q})_{+}$ does not have non-trivial normal
subgroups, which is known in some cases (see \cite[Ch. 9]{PR}).

iii) We need property c) only to distinguish between Shimura varieties and
equivariant fibrations over Shimura varieties, whose fibers are 
projective 
homogeneous spaces for the action of $\bo{G}_{\B{C}}$. 
Though we will not consider these objects in this paper, they are no
less interesting, than Shimura varieties themselves.

iv) We can replace property f), which is (one of) the strongest of 
our properties, by some weaker property f${}_0$), which will still
guarantee that the triple comes from the Shimura variety, corresponding to some
semi-simple adjoint group $\bo{G'}$, locally isomorphic to $\bo{G}$ at 
all places. For example, we can take f${}_0$) be the property f${}_1$)
from \rr{propf}, but a much weaker property will suffice. 

v) If we omit property f), then we can still guarantee that
for each sufficiently small compact and open subgroup 
$S\subset\bo{G}(\B{A}^{\be})$, each connected
component of the triple $(X_S;(P_S,\C{H}_S);\rho)$ is of the form
considered in \re{triple}.

b) Remark a) indicates that we can essentially ``classify'' all the triples which 
fail to satisfy some of the properties a), c), f) or g). On the other hand
we do not know how to  ``classify'' all the triples which 
fail to satisfy some of the properties b), d) and e). We assume that 
an answer to this question may lead to a better 
understanding of Shimura varieties.
\end{Rem}

As a first application of the Main Theorem we get the explicit form
of Kazhdan's theorem (or, equivalently, a weaker version of Langlands'
conjecture) on the conjugation of Shimura varieties (compare
\cite{Brv} and \cite{Mi1}).  
First we have to introduce some notation.

Fix an element $\sigma\in\Aut(\B{C})$. Then the conjugate
 ${}^{\sigma}\check{M}$ of $\check{M}$ is also a conjugate class of 
special parabolic subgroups of $\bo{G}_{\B{C}}$. Hence (see
 \re{dual} and \re{inn}) there exists a unique (up to 
a $\bo{G}(\B{C})$-conjugacy) pair consisting of an inner twist  
$\bo{\wt{G}}^{\sigma}$ of $\bo{G}_{\B{R}}$ and a 
$\bo{\wt{G}}^{\sigma}(\B{R})$-conjugacy class $M^{\sigma}$ of
 special homomorphisms $h:S^1\hra\bo{\wt{G}}^{\sigma}(\B{R})$ such that
 ${}^{\sigma}\check{M}$ coincides with the generalized Grassmanian 
corresponding to $(\bo{\wt{G}}^{\sigma},M^{\sigma})$.
Let $\wt{c}_{\sigma}\in H^1(\B{R},\bo{G})$ be the cohomology class,
corresponding to $\bo{\wt{G}}^{\sigma}$.

\begin{Cl} \label{C:coh}
There exists a unique cohomology class 
$c_{\sigma}\in H^1(\B{Q},\bo{G})$, whose image in $H^1(\B{R},\bo{G})$
is $\wt{c}_{\sigma}$ and whose image in $H^1(\B{Q}_p,\bo{G})$ is 
trivial for each finite prime $p$.
\end{Cl}
\begin{proof}
The uniqueness of $c_{\sigma}$ follows from the Hasse principle for
adjoint groups, thus it remains to show its existence.

By the proof of \rp{prop}, e), $\bo{G}$ contains a maximal elliptic
torus $\bo{T}$, defined over $\B{Q}$, and there exists an element $h\in M$ which factors 
through $\bo{T}(\B{R})$. Let $\mu$ be the corresponding cocharacter of 
$\bo{T}_{\B{C}}$. Since the element
$\sigma(\mu(-1))/\mu(-1)\in\bo{T}(\B{R})$ is of order two, it defines
a certain cohomology class $\wt{t}_{\sigma}\in H^1(\B{R},\bo{T})$.
Moreover, it follows from the construction that the image of $\wt{t}_{\sigma}$
in $H^1(\B{R},\bo{G})$ is $\wt{c}_{\sigma}$ (see \re{mu}). 
Hence it will suffice to show
the existence of a cohomology class
$t_{\sigma}\in H^1(\B{Q},\bo{G})$, whose image in $H^1(\B{R},\bo{T})$
is $\wt{t}_{\sigma}$ and whose image in $H^1(\B{Q}_p,\bo{T})$ is 
trivial for each finite prime $p$.

By the local and the global Tate-Nakayama duality (see \cite[$\S$3]{Ko}),
this would follow if we show
that $\wt{t}_{\sigma}$ belongs to the kernel of the natural composition map
$$\psi:H^1(\B{R},\bo{T})\isom H^1(\B{R},X^*(\bo{T}))^D\to
H^1(\B{Q},X^*(\bo{T}))^D,$$ 
where by $(\cdot)^D$ we denote the dual of a finite abelian group.

Since $\bo{T}_{\B{R}}$ is anisotropic, $\bo{T}$ splits over a certain 
Galois CM-extension $K/\B{Q}$. Then we have an action 
of the Galois group  $Gal(K/\B{Q})$ on $X^*(\bo{T})$.
Since $K$ is a CM-field, this action commutes with the complex
conjugation, and, therefore, induces an action on the 
cohomology group $H^1(\B{R},X^*(\bo{T}))$.
Using the definition of the cohomology groups, one can check that  
the map $\psi$ is constant on the orbits of  $Gal(K/\B{Q})$.

By the construction, the cohomology class $\wt{t}_{\sigma}$ is equal to
$\sigma(\wt{t})-\wt{t}$, where we denote
by $\wt{t}$ the cohomology class in $H^1(\B{R},\bo{T})\cong
H^1(\B{R},X^*(\bo{T}))^D$, induced by the
element $\mu(-1)\in\bo{T}(\B{R})$. Hence  
$\wt{t}_{\sigma}$ has to belong to the kernel of $\psi$, implying 
the assertion. 
\end{proof}

\begin{Not}
Let $\bo{G}^{\sigma}$ be the inner form of $\bo{G}$, corresponding to 
$c_{\sigma}$, and let $\Phi_{\sigma}$ be any isomorphism 
$\bo{G}_{\B{C}}\times\bo{G}(\B{A}^{\be})\isom
\bo{G}^{\sigma}_{\B{C}}\times\bo{G}^{\sigma}(\B{A}^{\be})$,
induced by the inner twisting $c_{\sigma}$. Notice that
$\Phi_{\sigma}$ is unique up to a conjugation.
\end{Not}

\begin{Cor} \label{C:conj}
Let $T=T(\bo{G}, M)$ be the triple, 
corresponding to a Shimura pair $(\bo{G}, M)$, and let 
$\sigma$ be an automorphism of $\B{C}$. Then the conjugate ${}^{\sigma}T$
 is  $\Phi_{\sigma}$-equivariantly 
isomorphic to the triple $T(\bo{G}^{\sigma}, M^{\sigma})$,
corresponding to the Shimura pair $(\bo{G}^{\sigma}, M^{\sigma})$.
\end{Cor}

\begin{pf}
 Conjugating the triple $T: (X; (P,\C{H});\rho:P\to\check{M})$ by $\sigma$,
 we get the triple 
${}^{\sigma}T:({}^{\sigma}X;
 ({}^{\sigma}P,{}^{\sigma}\C{H});{}^{\sigma}\rho:{}^{\sigma}P\to{}^{\sigma}\check{M})$,
 equipped with the natural action of the group
 ${}^{\sigma}\bo{G}_{\B{C}}\times\bo{G}(\B{A}^{\be})=\bo{G}_{\B{C}}\times\bo{G}(\B{A}^{\be})$.

Explicitly, $(g_{\B{C}},g_f)({}^{\sigma}p)={}^{\sigma}[({}^{\sigma^{-1}}g_{\B{C}}),g_f)(p)]$ for all
$g_{\B{C}}\in\bo{G}(\B{C})$, $g_f\in\bo{G}(\B{A}^{\be})$ and $p\in P(\B{C})$.
In particular, the stabilizer of ${}^{\sigma}p\in{}^{\sigma} P(\B{C})$ is
\begin{equation} \label{e:stab}
 \Gm_{{}^{\sigma}p}=\{({}^{\sigma}g_{\B{C}}, g_f)\in\bo{G}(\B{C})\times\bo{G}(\B{A}^{\be})\,|\,
(g_{\B{C}}, g_f)\in\Gm_p\}.
\end{equation}
By definition, the triple ${}^{\sigma}T$ is 
$\Phi_{\sigma}$-equivariantly isomorphic to a certain triple 
$T^{\sigma}$, consisting of data i)-iii) of \re{shimtr}
with respect to the Shimura pair $(\bo{G}^{\sigma}, M^{\sigma})$. 
By the Main Theorem it remains to check that $T^{\sigma}$
satisfies the properties a)-g) of \rp{prop}.

Using (\ref{e:stab}) we see that the properties a)-d) and g) 
are preserved under the conjugation and that e)
is preserved in the presence of f${}_1)$. Finally, since each
$\Gm_p\subset\bo{G}(\B{C})\times\bo{G}(\B{A}^{\be})$ is conjugate to a 
subgroup of $\bo{G}(\B{Q})$, embedded diagonally, the same is true for
each $ \Gm_{{}^{\sigma}p}$ by  (\ref{e:stab}). Since the isomorphism
$\Phi_{\sigma}$ is induced by an inner twisting, property f) for 
$T^{\sigma}$ follows as well.
\end{pf}

As a further application we show the existence of the canonical models
in the case of adjoint groups. 

\begin{Not}
For a Shimura pair  $(\bo{G}, M)$ 
we denote by $E(\bo{G}, M)\subset \B{C}$  
the field of rationality of the corresponding conjugacy class 
$\check{M}$ of special parabolic subgroups of $\bo{G}_{\B{C}}$, and we
will call it the {\em reflex field} of $(\bo{G}, M)$. 
\end{Not}

\begin{Cor} \label{C:can}
The Shimura variety  $X=X(\bo{G},M)$ and the full triple
$T=T(\bo{G},M)$ have 
unique $\bo{G_{E'}}\times\bo{G}(\B{A}^{\be})$-equivariant structures
over every subfield $E'$ of $\B{C}$, containing the reflex field  
$E=E(\bo{G}, M)$ of $(\bo{G}, M)$. Moreover, $E$ is the smallest
subfield of $\B{C}$ over which $T$ (resp. $X$) has an equivariant model.
\end{Cor}

\begin{Rem}
This is the only place in the paper where we use in a very strong way the assumption that $\bo{G}$ is 
adjoint.
\end{Rem}
\begin{pf} 
The definition of $E$ implies that for each $\sigma\in\Aut(\B{C}/E)$, the pair 
$(\bo{G}^{\sigma}, M^{\sigma})$ is conjugate to
$(\bo{G}, M)$ inside $\bo{G}_{\B{C}}$. Hence each  isomorphism
$\Phi_{\sigma}$ can be chosen to be the identity map.
Thus  \rco{conj} implies that there exists 
$\bo{G}_{\B{C}}\times\bo{G}(\B{A}^{\be})$-equivariant isomorphisms
$\varphi_{\sigma}:{}^{\sigma}T\isom T$ and 
$\phi_{\sigma}:{}^{\sigma}X\isom X$.
Moreover, since $\bo{G}$ is adjoint, both $X$ and $T$ have no
equivariant isomorphisms (see \rl{aut},2)-4)), hence
each $\varphi_{\sigma}$ and $\phi_{\sigma}$ is unique. 
As a consequence, the equivariant structures of $T$ and $X$ are
unique, and the $\varphi_{\sigma}$'s 
satisfy the cocycle condition, thus forming the descent datum.

To show that this descent datum is effective (see Appendix) we observe 
that the $\varphi_{\sigma}$'s define descent datum on each 
$X_S$. Moreover, these descent data commute with the natural projections
$\pi_{S,T}:X_S\to X_T$ for all $S\subset T$ and with isomorphisms
$g_S:X_S\isom X_{gSg^{-1}}$ for all $g\in\bo{G}(\B{A}^{\be})$. 
Applying now \rco{aut} to the system $\{X_S\}_S$ with morphisms
$\{\pi_{S,T}\}_{S,T}$ and  $\{g_S\}_{g,S}$ we conclude from 
\rl{aut}, 1),2) that $X$ descends effectively to $E$. 
Knowing the effectivity of the descent datum for $X$ we conclude 
from \rp{cont} and from 3),4) of \rl{aut} that the whole triple 
equivariantly descends to $E$.

Conversely, assume that ${}^{\sigma}X$ is 
$\bo{G}(\B{A}^{\be})$-equivariantly isomorphic to $X$ for some 
$\sigma\in\Aut(\B{C})$. \rco{conj} then implies that 
$X(\bo{G}^{\sigma},M^{\sigma})$ is $\Phi_{\sigma}$-equivariantly
isomorphic to $X(\bo{G},M)$. By the proof of \rr{isom}, b),
some conjugate of $\Phi_{\sigma}$ is induced by a certain isomorphism
$\bo{G}\isom\bo{G}^{\sigma}$, which maps $M$ into $M^{\sigma}$.
As $\Phi_{\sigma}$ is induced by an inner twisting, we get that the pairs
$(\bo{G}^{\sigma}_{\B{R}}, M^{\sigma})$ and 
$(\bo{G}_{\B{R}}, M)$ are conjugate inside $\bo{G}_{\B{C}}$.
Therefore ${}^{\sigma}\check{M}=\check{M}$, hence 
$\sigma$ belongs to $\Aut(\B{C}/E)$, implying the minimality property of $E$.
\end{pf}

%

a further less obvious application of the Main Theorem is the following theorem on 
the $p$-adic uniformization (proven in \cite{Va1,Va2}, where the reader
is reffered  for a more precise and more general statement and 
for a detailed proof).

\begin{Cor} \label{C:p-adic}
Let $\bo{G}/\B{Q}$ be a group of unitary similitudes modulo center
such that $\bo{G}(\B{R})\cong\bo{PGU_{d-1,1}}(\B{R})^r\times\bo{PGU_d}(\B{R})^{g-r}$ and 
$\bo{G}(\B{Q}_p)\cong\prod_{i=1}^r PD\m_{L_i,1/d}\times \wt{G}_p$, where $\wt{G}_p$ is a locally compact 
topological group, and for each $i=1,...,r$ the field $L_i$ is a
finite extension of $\B{Q}_p$ and $D_{L_i,1/d}$ is a division algebra 
over $L_i$ with Brauer invariant $1/d$.  

Let  $E$ be the reflex field the Shimura pair 
$(\bo{G}, M)$ for a suitable $M$, and let 
$v$ be a prime of $E$, which lies over $p$ and satisfies some natural 
conditions (see \cite[2.5]{Va2}). In particular,
we assume that $E_v$ is the composite field of the images of the
$L_i$'s under certain embeddings $L_i\hra\B{C}_p$. 

Let $X$ be the unique
$\bo{G}(\B{A}^{\be})$-equivariant structure over $E$ (see \rco{can})
of the complex Shimura variety $X(\bo{G}, M)$, and
let $S\subset\bo{G}(\B{A}^{\be})$ be an open and compact subgroup of the form 
$\prod_{i=1}^r PD\m_{L_i,1/d}\times S'$ for some 
$S'\subset\wt{G}_p\times\bo{G}(\B{A}^{\be,p})$. Then
\begin{equation} \label{E:pad}
(S\bs X\otimes_E E_v)\an\cong S'\bs[\prod_{i=1}^r \Omega^d_{L_i}\otimes_{L_i} E_v\times 
(\wt{G}_p\times\bo{G}(\B{A}^{\be,p})]/\bo{G'}(\B{Q}),
\end{equation}
where 

a) $\Omega^d_{L_i}$ is the Drinfeld $p$-adic upper half space over
$L_i$ (an open analytic subset of $\B{P}^{d-1}_{L_i}$, obtained by 
removing from it the union of all $L_i$-rational hyperplanes).

b) $\bo{G'}$ is an inner form of $\bo{G}$ such that
$\bo{G'}(\B{A}^{\be,p})\cong\bo{G}(\B{A}^{\be,p}),\;\bo{G'}(\B{R})\cong\bo{PGU_d}(\B{R})^g$ and
$\bo{G'}(\B{Q_p})\cong\prod_{i=1}^r\bo{PGL_d}(L_i)\times \wt{G}_p$.
In particular, 
$\bo{G'}(\B{A}^{\be})\cong\prod_{i=1}^r\bo{PGL_d}(L_i)\times\wt{G}_p\times\bo{G}(\B{A}^{\be,p})$. 
\end{Cor}

\begin{Ex}
(very particular case of \rco{p-adic}) Let $B$ be a quaternion algebra over $\B{Q}$, splitting at $\be$ and ramifying at $p$.
Let $\bo{G}=\bo{PGL_1}(B)$, and let $X$ be the Shimura curve 
corresponding to $\bo{G}$. Then the corresponding reflex field is
$\B{Q}$, and $(S\bs X\otimes_{\B{Q}}\B{C})\an\cong S\bs
[(\B{P}^1(\B{C})\sm\B{P}^1(\B{R}))\times \bo{G}(\B{A}^{\be})]/\bo{G}(\B{Q})$
is a compact curve. \rco{p-adic} asserts in this case that 
$$(S\bs X\otimes_{\B{Q}}\B{Q}_p)\an\cong S'\bs[\Omega^2_{\B{Q}_p}\times\bo{G}(\B{A}^{\be,p})]/P(B')\m,$$
where $B'$ is the unique quaternion algebra over $\B{Q}$, which splits at $p$, ramifies at $\be$ and is
locally isomorphic to $B$ at all the other places of $\B{Q}$.  
\end{Ex}

\begin{pf}
(Sketch) Recall first that for every cocompact lattice 
$\Gm\subset\prod_{i=1}^r\bo{PGL_d}(L_i)$ the quotient
$\wt{Y}_{\Gm}:=\Gm\bs \prod_{i=1}^r \Omega^d_{L_i}\otimes_{L_i} E_v$ 
has a unique structure of a 
projective scheme over $E_v$. As the right hand side of (\ref{E:pad}) 
is a finite disjoin union of analytic spaces of the form
$\wt{Y}_{\Gm}$, it also has a unique structure $Y_{S}$ of a 
projective scheme over $E_v$.
As in the case of Shimura varieties, we can form an inverse limit $Y'$
of the $Y_{S}$'s, equipped
with an action of $\wt{G}_p\times\bo{G}(\B{A}^{\be,p})$. Moreover,
using Drinfeld's coverings 
$\Sigma^d_{L_i}$'s of $\Omega^d_{L_i}$'s we can construct a  
$\wt{G}_p\times\bo{G}(\B{A}^{\be,p})$-equivariant Galois cover $Y$ of
$Y'$ with a Galois group $\prod_{i=1}^r PD\m_{L_i,1/d}$. As a consequence,
$Y$ is equipped with an action of the group $\bo{G}(\B{A}^{\be})$, and 
$S\bs Y\cong Y_S$ for each $S$ as in the corollary. 
Thus it is enough
to show that $Y$ is $\bo{G}(\B{A}^{\be})$-equivariantly isomorphic to 
$X_{E_v}$.

As in the complex case, we construct a triple over $E_v$ 
consisting of $Y$, a $\bo{G}_{E_v}$-torsor $P$ on $Y$ equipped with an 
equivariant flat connection on it and an equivariant morphism from $P$
to $(\B{P}_{E_v}^{d-1})^r$. After extending scalars to $\B{C}$ 
(using an arbitrary field embedding of $E_v$ into $\B{C}$), 
we get a triple, satisfying assumptions i)-iii) of \re{shimtr}. 

If we show that this triple satisfies all the properties a)-g) of 
\rp{prop}, the Main Theorem would imply that 
$Y$ and $X_{E_v}$ are $\bo{G}(\B{A}^{\be})$-equivariantly isomorphic over
$\B{C}$. By \rc{can} they would be then actually isomorphic over
$E_v$, implying the assertion. The proof of the properties a)-g)
is very similar to that of \rp{prop}, and is given in detail in 
\cite{Va1,Va2}.  
\end{pf}

\section{proof of the Main Theorem} \label{S:proof}

The strategy of the proof is very similar to that described in
\cite{Va1, Va2}. For the convenience of the reader we describe all the
steps in the proof but will omit some of the technical detailes for which we will reffer
to the corresponding places in these papers.

{\bf Step 0}. Reduction of the problem.

I) Since the transition maps $X_T\to X_S$ are finite and etale for all 
$T\subset S$ sufficiently small (by property \re{shimtr}, i)), there exists an inverse limit $X\an$ of 
the $X_S\an$'s in the category of complex analytic spaces 
(compare \cite[Lem. 1.6.3]{Va1}).

II) Let $M'$ be a connected component of $X\an$, and let 
$\Dt\subset\bo{G}(\B{A}^{\be})$ be the stabilizer of $M'$. Then property
f${}_4$) (see \rr{propf}) implies that $\Dt$ acts faithfully on $M'$.
Since for each complex algebraic variety $Y$ the set of the connected 
components of $Y$ coinsides with that of $Y\an$, property a) 
implies that the group $\bo{G}(\B{A}^{\be})$ acts transitively on the set of the
connected components of $X\an$ (compare \cite[the proof of
Prop. 1.6.1]{Va1}), so that 
$X\an\cong[M'\times\bo{G}(\B{A}^{\be})^{disc}]/\Dt$.

III) Let $\pi:\wt{M}\to M'$ be the universal cover, and let 
$\wt{\Dt}\subset\Aut(\wt{M'})$ be the group of all liftings of elements
of $\Dt\subset\Aut(M)$. Then $X\an\cong[\wt{M}\times\bo{G}(\B{A}^{\be})^{disc}]/\wt{\Dt}$,
and  $P\an\cong[\wt{Y}\times\bo{G}(\B{A}^{\be})^{disc}]/\wt{\Dt}$ for
a certain $(\bo{G}_{\B{C}})\an$-principal bundle $\wt{Y}$
over $\wt{M}$, equipped with a $\wt{\Dt}$-invariant flat connection.    
Since $\wt{M}$ is simply connected, $\wt{Y}$ naturally
decomposes as a product $\wt{M}\times(\bo{G}_{\B{C}})\an$, and
we get a $\bo{G}_{\B{C}}\times\bo{G}(\B{A}^{\be})$-equivariant isomorphism 
 $$P\an\cong[\wt{M}\times\bo{G}_{\B{C}}\an\times\bo{G}(\B{A}^{\be})^{disc}]/\wt{\Dt},$$
where $\wt{\Dt}$ acts on the second factor via some homomorphism
$j:\wt{\Dt}\to\bo{G}(\B{C})$ and  where $\C{H}\an$ corresponds to the 
trivial connection on $\wt{M}\times\bo{G}_{\B{C}}\an$ 
(compare \cite[Prop. 4.4.2]{Va1}). Moreoevr, $j$ is unique up to a conjugacy
 

V) Let $\rho_0$ be the restriction of $\rho\an:P\an\to\check{M}\an$
to $\wt{M}\cong\wt{M}\times\{1\}\times\{1\}$. Then 
$\rho_0:\wt{M}\to\check{M}\an$ is a $\wt{\Dt}$-equivariant map, which
is a local isomorphism by property \re{shimtr}, iii) of the triple.

VI) By statements II)-IV) above, we have a natural embedding of 
$\wt{\Dt}$ into $\Aut(\wt{M})\times\bo{G}(\B{C})\times\bo{G}(\B{A}^{\be})$.
Therefore the Main Theorem asserts that

 $(*)$ both $M'$ and $\wt{M}$ are  isomorphic to a connected 
component of $M$, and $\rho_0$ corresponds under this isomorpism 
to the Harish-Chandra-Borel embedding, and 

 $(**)$ we have isomorphisms $\wt{\Dt}\isom\Dt\isom\bo{G}(\B{Q})_+$, 
and our embedding of $\wt{\Dt}$ into 
$\bo{G}(\B{C})\times\bo{G}(\B{A}^{\be})$ corresponds under this 
isomorphism to the diagonal embedding of $\bo{G}(\B{Q})_+$ 
(up to a conjugation).

Moreover, it follows from the uniqueness of the algebraic structure of 
our triple that the Main Theorem actually follows from $(*)$ and $(**)$.

{\bf Step 1}. Proof of the statement $(*)$ of Step $0$ (compare \cite[Section 3,
Step 2]{Va2} or \cite{NR}).

I) Let $\wt{J}$ be the closure of $\wt{\Dt}$ in $\Aut(\wt{M})$, 
let $j_0:\wt{\Dt}\to\Aut(\check{M}\an)$ be the composition of $j$ with 
the natural map $\bo{G}(\B{C})\to\Aut(\check{M}\an)$,
and let $J$ be the closure of $j_0(\wt{\Dt})$. Using the fact that 
$\rho_0$ is a local homomorphism, we get that $j_0$ extends by
continuity to a continuous homomorphism $\phi:\wt{J}\to J$ 
(see \cite[Lem. 3.1]{Va2}).

II) Let $\bo{G_0}$ be the kernel of the action of $\bo{G}$ on 
$\check{M}$, then $\bo{\bar{G}}:=\bo{G}/\bo{G_0}$ is isomorphic to 
$\Aut(\check{M})^0$. Property b) implies that $j(\wt{\Dt})$ is Zariski 
dense in $\bo{G}$, therefore $J$ is a Lie subgroup of 
$\bo{\bar{G}}(\B{C})$, which is Zariski dense in $\bo{\bar{G}}$. 

III) Let $\bo{G_1}, \bo{G_2},...,\bo{G_r}$ be all the simple factors of 
$\bo{\bar{G}}$. Properties d) and e) imply that $J$ contains a maximal 
compact torus $T$ of  $\bo{\bar{G}}(\B{C})$, stabilizing 
$\rho_0(y)\in\check{M}$ for some point $y\in\wt{M}$.
(Actually the above properties imply that $T$ is contained in 
the closure of $j(Stab_{\wt{\Dt}}(y))$).
 This together with II) imply by standard argument involving
Lie algebras (see \cite[3.3]{Va2}) that 
$J^0=\prod_{i=1}^r J_i\subset\prod_{i=1}^r \bo{G_i}(\B{C})$, where each
$J_i$ is either $\bo{G_i}(\B{C})$ or $\bo{J_i}(\B{R})^0$ for one of 
the real forms $\bo{J_i}$ of $\bo{G_i}$. In both cases 
$J^0\cong\bo{J}(\B{R})^0$ for some semi-simple adjoint real group $\bo{J}$.

IV) Choose a $T$-invariant Riemannian metric on $\check{M}\an$ (for example,
the one making it a Hermitian symmetric space). Then its pull-back 
under $\rho_0$ is a Riemannian metric on $\wt{M}$,
invariant under under $\phi^{-1}(T)$. Hence the intersection 
$\phi^{-1}(T)\cap Stab_{\wt{J}}(y)$ is a closed subgroup
of the isotropy group at $y$, therefore it is compact. 
It follows that the closure of $Stab_{\wt{\Dt}}(y)$ in $\wt{J}$ is
compact as well. Our choice of $T$ and $y$ therefore implies that
the image of  $\overline{Stab_{\wt{J}}(y)}$ under $\phi$ contains $T$. 
Since $\Lie\, T$ spans all of $\Lie\,J$ as an $Ad j_0(\wt{\Dt})$-module
over $\B{R}$ and since $\phi$ has a dense image, $\phi$ has to be 
surjective. Since $\rho_0$ is a local isomorphism, we get that 
$\wt{J}$ is a Lie group, and that $\phi$ is a topological covering.

V) Using again property d) together with some basic structure theorems
of reductive algebraic groups we show (as in \cite[Prop. 3.5]{Va2}) 
that the stabilizer $J_y$ of $\rho_0(y)\in\rho_0(\wt{M})$ is compact 
for each point $y\in\wt{M}$ and that $\bo{J}$ is a real form of 
$\overline{\bo{G}}_{\B{C}}$. In particular, each
$J_y$ is contained in some maximal compact subgroup 
$ J^{\text{comp}}\subset\overline{\bo{G}}(\B{C})$.

VI) It now formally follows from V) by dimension counting that the group
$J$ is connected, that it acts transitively on $\rho_0(\wt{M})$,
that $\wt{J}$ acts transitively on $\wt{M}$, and that 
$\rho_0:\wt{M}\to\rho_0(\wt{M})$ is a topological covering 
(compare \cite{Va2}).

VII) By the characterization of Hermitian symmetric spaces we see
that $\rho_0(\wt{M})$ is a Hermitian symmetric space, 
not containing Euclidean factors. Moreover, $\rho_0(\wt{M})$
can be decomposed as a product $\prod_{i=1}^r M'_i\subset
\prod_{i=1}^r \check{M_i}\an$, where  
$\check{M}\an=\prod_{i=1}^r \check{M_i}\an$ is the 
decomposition of $\check{M}\an$ into the product of irreducible 
Hermitian spaces and each embedding $M'_i\hra\check{M_i}\an$
is either a Harish-Chandra-Borel embedding or an isomorphism.
Moreover, property c) implies that $\rho_0(\wt{M})$ can not have compact
factors, therefore is has to be isomorphic to a connected component
$M^0$ of $M$.

VIII) Since $M^0$ is simply connected, we get
$\wt{M}\cong\rho_0(\wt{M})\cong M^0$. As the group $\Aut(M^0)$ has a trivial 
center, it has no non-trivial normal discrete subgroups. This implies 
that $M'\cong M$ (compare \cite[the proof of Prop. 1.6.1]{Va1}),
hence $\wt{\Dt}\cong\Dt$.
This completes the proof of the statement $(*)$ of Step $0$.

{\bf Step 2}. Proof of the statement $(**)$ of Step $0$.

I) Recall first that $\Dt$ is naturally embedded into 
$\Aut(M^0)^0\times\bo{G}(\B{C})\times\bo{G}(\B{A}^{\be})$ and 
that its projection to $J\cong\Aut(M^0)^0$ coincides with $j_0$.
We claim that the projection of $\Dt$ to the product of the first and
the third factors is a lattice in $J\times\bo{G}(\B{A}^{\be})$. 

It will suffice to show that $\Dt_S:=\Dt\cap S$ is a lattice in $J$
for some sufficiently small open and compact subgroup 
$S\subset\bo{G}(\B{A}^{\be})$ (compare \cite[the proof of
Prop. 1.6.1]{Va1}). By Step $1$, $M^0$ is the universal cover of the 
 the quotient $Y_S:=\Dt_S\bs M^0$, hence $\Dt_S$ is discrete in   
$\Aut(M^0)^0\cong J$.
Since the stabilizer of every point in $M^0$ is
compact, it remains to show that $Y_S=\Dt_S\bs M^0$ 
has a finite volume with respect to some (or, equivalently, to every) 
$J^0$-invariant non-trivial measure on $M^0$. 
The statement follows formally from
the fact that $Y_S$ is a complex analytic space, corresponding to a 
complex quasi-projective variety (a connected 
component of $X_S$) and that Kobayashi hyperbolic measure on $M^0$ is 
non-trivial (compare \cite{Ka2}).

Indeed, by Hironaka's theorem on the resolution of singularities, 
$Y_S$ can be embedded into a compact complex manifold $\overline{Y}_S$ 
such that $D:=\overline{Y}_S\sm Y_S$ is a divisor with normal crossings.
Then for every $y\in\overline{Y}_S$ there exists an open neighborhood 
$U_y\subset\overline{Y}_S$ of $y$ and a biholomorphic map 
$\varphi_y$ between $U_y$ and the unit polydisc, 
which maps $y$ into the center of the polydisc such that
$\varphi_y$ maps $U'_y:=U_y\cap Y_S$ into a product of unit discs and 
punctured unit discs. Let $V_y\subset U_y$ be the preimage under
$\varphi_y$ of the product of open
unit discs of radius $1/2$, and set $V'_y:=V_y\cap U'_y$.
Since $\overline{Y}_S$ is compact, there exits a finitely many points 
$y_1,...,y_k\in\overline{Y}_S$ such that  $\overline{Y}_S$ is covered by
the $V_{y_i}$'s, hence $Y_S$ is covered by the  $V'_{y_i}$'s.

For every complex manifold $N$ we denote by $\mu_N$ the Kobayashi hyperbolic
measure on $N$ (see \cite[IV, $\S$4]{La}). Since $M^0$ is isomorphic
to a bounded symmetric domain, $\mu_{M^0}$ is a non-trivial  
$J^0$-invariant measure on $M^0$. 
Moreover, by the properties of the Kobayashi measures we have

$$\mu_{M^0}(Y_S)=\mu_{Y_S}(Y_S)\leq\sum_{i=1}^k\mu_{Y_S}(V'_{y_i})
\leq\sum_{i=1}^k\mu_{U'_{y_i}}(V'_{y_i}),$$
and direct calculation shows that each summand of the last expression is 
finite (compare \cite[II, Prop 3.1]{La}). This shows the assertion.

II) For each finite set $N$ of non-archimedean primes 
of $F$ and for each  compact and open subgroup 
$S\subset\bo{H}(\B{A}_F^{\be,N})$ we consider a lattice $\Dt^{N;S}
\subset\bo{J}(\B{R})\times\prod_{u\in N}\bo{H}(F_u)$,  
consisting of projections of elements of $\Dt$, whose components outside 
$N$ belong to $S$. Let $N_{\text{is}}$ (resp. $N_{\text{unis}}$) 
be the set of all $u\in N$ such that $\bo{H}(F_u)$ is non-compact 
(resp.  $\bo{H}(F_u)$ is compact). 

III) By the same arguments as in \cite[Lem. 3.9]{Va2} and 
using the fact that the set of regular elliptic elements of $J$ is open
(compare the proof of \rp{prop}, d)), we deduce from condition
f${}_1$) of \rr{propf}  
that $\Dt$ projects injectively on each simple factor 
of $\bo{J}(\B{R})\times\bo{H}(\B{A}_F^{\be})$. In particular, 
the lattice $\Dt^{N;S}$ is irreducible. Furthermore, 
 $\Dt^{N;S}$ has an infinite index in its commensurator, since it has 
an infinite index in $\Dt$. Therefore \cite[Thm. 1.13, p.299]{Ma} 
implies that the projection of $\Dt^{N;S}$ to 
$\bo{J}(\B{R})\times\prod_{u\in N_{\text{is}}}\bo{H}(F_u)$ is arithmetic.

IV) In the case that $N_{\text{unis}}$ is non-empty,  
we can not conclude directly that  $\Dt^{N;S}$ is arithmetic, since we do 
not know yet that it satisfies property (QD) of \cite[p.289]{Ma}.
In this case we will argue as follows:
Fix   $v\in N_{\text{unis}}$. Let $\bo{H_v}\subset\bo{H_{F_v}}$ 
be the connected component of the identity of the Zariski closure of
the projection of  
$\Dt^{N;S}$ to $\bo{H}(F_v)$, and let $\Dt'_v\subset\Dt^{N;S}$ be the 
inverse image of  $\bo{H_v}$. Then $\Dt'_v$ has a finite index in 
$\Dt^{N;S}$. Since $\bo{H}(F_v)$ is a compact group, 
it does not contain unipotent elements over $F_v$, hence $\bo{H_v}$ is reductive,. 
Moreover, $\bo{H_v}$ is a non-trivial centerless group,
because it is a Zariski closure of a centerless group $\Dt'_v$. Then the 
standard Lie algebra argument shows that the closure of $\Dt'_v$ in  
$\bo{H_v}(F_v)$ is open.
\cite[Thm. 1.13, p.299]{Ma} now implies that the projection of 
$\Dt'_v$ to $\bo{J}(\B{R})\times\prod_{u\in N_{\text{is}}}\bo{H}(F_u)\times
\bo{H_v}(F_v)$ is arithmetic. In particular, $\bo{H_v}$ is a form of
$\bo{H}$, hence $\bo{H_v}=\bo{H_{F_v}}$ and $\Dt'_v=\Dt^{N;S}$.

V) By the definition of arithmetic groups, 
there exists a number field $F'$, an absolutely simple group 
$\bo{H'}$ over $F'$, a map $I=I_N$ from $N\cup\cup_{i=1}^r \be_i$ to the set 
of places of $F'$ and a continuous isomorphism
$$\varphi_N:\bo{J}(\B{R})\times\prod_{u\in N}\bo{H}(F_u)\isom
\prod_{i=1}^r\bo{H'}(F'_{I(\be_i)})\times\prod_{u\in N} \bo{H'}(F'_{I(u)}),$$
induced by the product of local algebraic isomorphisms (compare 
for example \cite[2.4.6]{Va1} for a more precise statement) such that
$\varphi_N(\Dt^{N;S})\subset\bo{H'}(\C{O}_{F'})$  is commensurable 
with $\bo{H'}(\C{O}_{F';I(N)})$,
embedded diagonally. (Here by $\C{O}_{F';I(N)}$ we denote the set of
elements of $F'$, which are integral outside $I(N)$.) 
Moreover, $F'$ and $\bo{H'}$ do not depend on $N$ and $S$. 

Taking $N$ larger and larger, we find a (non-necessary injective) map 
$I$ from the union of $\cup_{i=1}^r\be_i$ and the set 
of all non-archimedean places of $F$ to the set of places of $F'$
and a continuous isomorphism 
$$\varphi:\bo{J}(\B{R})\times\prod_u\bo{H}(F_v)\isom
\prod_{i=1}^r\bo{H'}(F'_{I(\be_i)})\times\prod_u \bo{H'}(F'_{I(u)})$$
such that $\varphi(\Dt)\subset\bo{H'}(F')$. 

Finally notice that the definition of the arithmetic subgroup also
implies that $\bo{H'}$ is anisotropic at all the archimedean places of
$F'$, which does not belong to the image of $I$. 

VI) Next we claim that $F'$ is isomorphic to $F$. For example, it can be
seen as follows. \cite[Ch. VIII, Prop. 3.22]{Ma} shows that
$F'$ is equal to the field $\B{Q}(\TrAd\Dt_{\be_1})$, generated by
traces of the adjoint representation of  $\bo{J_1}$,
where by $\Dt_{\be_1}$ we denote the projection of $\Dt$ to 
$\bo{J_1}(\B{R})$. Using the fact that we have sufficiently many
regular elliptic elements in $\Dt\subset J$ (compare the proof of 
\rp{prop}, d)) we get, as in \cite[Lem. 2.3.3]{Va1}, that
\begin{equation} \label{E:TrAd}  
F'=\B{Q}(\TrAd\dt_{\be_1}|\dt\in\Dt\text{ is regular elliptic }).
\end{equation}
Property f${}_1$) then implies that $F'\subset F$. On the other hand,
we know by V), that $F\otimes_{\B{Q}}\B{Q}_p$ is a direct factor
of $F'\otimes_{\B{Q}}\B{Q}_p$ for almost all $p$. Hence $F'$ coincide
with $F$. Moreover, using again property  f${}_1$) we conclude from 
(\ref{E:TrAd}) and the definition
of the map $I$ (as in \cite[2.5.2]{Va1}) that $I$ is the identity map. 
(In particular, $I$ is one-to-one, showing a posteriori that
each $\Dt^{N;S}$ is arithmetic (compare the beginning of IV)).


VII) Next we show, following the arguments of \cite[Prop. 2.5.4]{Va1}, 
that $\varphi$ induces a continuous isomorphism 
$\bo{H}(\B{A}_F^{\be})\isom\bo{H'}(\B{A}_F^{\be})$.
Indeed, let $N_0$ be a finite set of finite primes of $F$ such that
both $\bo{H}$ and  $\bo{H'}$ are quasi-split outside $N_0$.
It remains to show that  $\varphi$ induces a continuous isomorphism 
$\bo{H}(\B{A}_F^{\be,N_0})\isom\bo{H'}(\B{A}_F^{\be,N_0})$. Since each lattice
$\Dt^{N;S}$ is irreducible, we conclude using \cite[Ch. II, Thm. 6.7 and 
Ch. I, Thm. 1.5.6]{Ma} that the closure of the projection of $\Dt$ to 
$\bo{H}(\B{A}_F^{\be,N_0})$ contains $\bo{H}(\B{A}_F^{\be,N_0})^{der}$  
(this is the analog of the strong approximation theorem in our situation). 
Then (compare \cite{Va1}) $\varphi$ extends by continuity to a continuous
isomorphism  
$\bo{H}(\B{A}_F^{\be,N_0})^{der}\isom\bo{H'}(\B{A}_F^{\be,N_0})^{der}$. 

The statement will follow if we show that the normalizer of 
$\bo{H}(\B{A}_F^{\be,N_0})^{der}$ inside $\prod_{u\notin N_0}\bo{H}(F_u)$
equals $\bo{H}(\B{A}_F^{\be,N_0})$ (and similarly for $\bo{H'}$), or, 
equivalently, that for each finite place
$u\notin N_0$ the normalizer of $\bo{H}(\C{O}_{F_u})^{der}$ inside 
$\bo{H}({F_u})$ equals $\bo{H}(\C{O}_{F_u})$.
But the compact group $\bo{H}(\C{O}_{F_u})^{der}$ has a
trivial centralizer inside $\bo{H}({F_u})$,
therefore its normalizer is a compact subgroup of $\bo{H}({F_u})$.
Since the normalizer clearly contains $\bo{H}(\C{O}_{F_u})$, which is
a maximal compact subgroup of $\bo{H}({F_u})$, it has to be equal to
it, as was claimed.

VIII) Put $\bo{G'}=\bo{R}_{F/\B{Q}}\bo{H'}$. 
By I), $\Dt$ is a lattice in $J\times\bo{G}(\B{A}^{\be})$. Hence V)
and VII) imply that $\varphi(\Dt)$ has a finite covolume in
$\bo{G'}(\B{A})$. Since $\bo{G'}(\B{Q})_{+}\subset\bo{G'}(\B{A})$ is
discrete, $\varphi(\Dt)$ has a finite index in  
$\bo{G'}(\B{Q})_{+}$. We claim that property g) implies that 
  $\varphi(\Dt)=\bo{G'}(\B{Q})_{+}$ (compare \cite[3.12]{Va2}).

Fix a regular elliptic element $\dt\in\bo{G'}(\B{Q})_{+}\subset J$
(compare the proof of \rp{prop}, d).
Its centralizer $\bo{T_{\dt}}$ is a maximal torus in $\bo{G'}$ such that
$\dt\in\bo{T_{\dt}}(\B{Q})$. Since $\bo{T_{\dt}}(\B{R})\subset J$ 
stabilizes a point in $M$, the image of
$\varphi(\Dt)\cap\bo{T_{\dt}}(\B{Q})$ in
$\bo{G}(\B{C})\times\bo{G}(\B{A}^{\be})$ stabilizes a point in $P$.
As $\varphi(\Dt)\cap\bo{T_{\dt}}(\B{Q})$ has a finite index in
$\bo{T_{\dt}}(\B{Q})$, it is Zariski dense in 
$\bo{T_{\dt}}$.
Hence property g) implies that $\bo{T_{\dt}}(\B{Q})\subset\varphi(\Dt)$.
Therefore $\dt$ belongs to $\varphi(\Dt)$. But $\dt$ was chosen
arbitrary, hence all regular elliptic elements 
$\dt\in\bo{G'}(\B{Q})_{+}\subset J$ belong to  
$\varphi(\Dt)$. But such elements generates 
$\bo{G'}(\B{Q})_{+}$ as a group, since regular elliptic elements
generate $J$ as a group, and since  
$\bo{G'}(\B{Q})_{+}$ is dense in $J$ (compare \cite[IX, Lem. 3.3]{Ma}). 
This shows the assertion.

IX) Recall that property b) implies that the image of the homomorphism 
$j:\bo{G'}(\B{Q})_{+}\to\bo{G}(\B{C})$ is Zariski dense, therefore the same
arguments as in \rr{propg} imply (using property f${}_2)$)
that $j$ extends to an algebraic isomorphism 
$\bo{G'}_{\B{C}}\isom\bo{G}_{\B{C}}$. Moreover, it follows from
the fact that  $\rho_0$ is a Harish-Chandra-Borel embedding and from the
definitions of $\bo{H}$ and $\bo{H'}$ that $j$ induces isomorpisms
$\bo{H}(F_{\be_i})^0\isom J_i^0\isom\bo{H'}(F_{\be_i})^0$ for each
$i=1,...,r$. Moreover, using V) and VI) 
we see that for all the archimedean places $\be'$ of
$F$, different from $\be_1,...,\be_r$, both $\bo{H}(F_{\be'})$  and 
$\bo{H'}(F_{\be'})$ are maximal compact subgroups of 
$\bo{H}(\B{C})\cong\bo{H'}(\B{C})$.
Since maximal compact subgroups of a complex reductive group
are all conjugate, we may therefore assume after possible replacing $j$ 
by its conjugate that $j$ induces an isomorphism 
$\bo{G'}(\B{R})\isom\bo{G}(\B{R})$.

X) We have proven in VII)-IX) the existence of a continuous isomorphism 
$\varphi:\bo{G}(\B{A})\isom\bo{G'}(\B{A})$ which is the product
of the local isomorphisms 
$\varphi_p:\bo{G}(\B{Q}_p)\isom\bo{G'}(\B{Q}_p)$ for all primes $p\leq\be$
(with $\B{Q}_{\be}:=\B{R}$) such that $\varphi(\Dt)=\bo{G'}(\B{Q})_{+}$.
Hence to finish the proof of the Main Theorem it remains to show that 
some conjugate of $\varphi$ is induced by a certain $\B{Q}$-rational isomorphism 
$\bo{G}\isom\bo{G'}$.

Choose first any isomorphism 
$\psi:\bo{G'}_{\overline{\B{Q}}}\isom\bo{G}_{\overline{\B{Q}}}$.
For each prime $p\leq\be$ choose an embedding 
$\overline{\B{Q}}\hra\overline{\B{Q}_p}$, and let
$\psi_p:\bo{G}_{\overline{\B{Q}_p}}\isom\bo{G}_{\overline{\B{Q}_p}}$
be the composition map $\psi\circ\varphi_p$. Replacing $\psi$ 
if necessary, we may 
assume that $\psi_{\be}$ is an inner automorphism.

We claim that each $\psi_p$ is an inner automorphism. This would
follow if we show
that for each conjugation-invariant regular function $f$ on $\bo{G}$ 
we have $f(\psi_p(g))=f(g)$ for each $g\in\bo{G}$. 
Choose first $g$ be the projection to $\bo{G}(\B{Q}_p)$ of an element 
$\dt\in\Dt$, whose projection to $J$ is elliptic. As in VIII), the image
of $\dt$ to $\bo{G}(\B{A}^{\be})\times \bo{G}(\B{C})$ stabilizes some
point in $P$.  
Then property f) implies that
$f(g)=f(\varphi^{-1}_{\be}\circ\varphi_p(g))$,
and the latter expression is equal to
$f(\psi^{-1}_{\be}\circ\psi_p(g))=f(\psi_p(g))$, 
since $\psi_{\be}$ is an inner automorphism. 
Since the projections of elliptic elements are
Zariski dense, the required equality holds for all $g\in\bo{G}$.

Now we are ready to prove our statement. Consider the cocycle
$\{\phi_{\sigma}\}_{\sigma}$ of $\Gal(\overline{\B{Q}}/\B{Q})$
in $\Aut(\bo{G})$, defined by
$\phi_{\sigma}:={}^{\sigma}\psi\circ\psi^{-1}$. We want to show
that it defines the trivial inner twisting. 
For each $p\leq\be$ and each 
$\sigma\in\Gal(\overline{\B{Q}_p}/\B{Q}_p)\subset\Gal(\overline{\B{Q}}/\B{Q})$
we have $\phi_{\sigma}={}^{\sigma}\psi_p\circ\psi_p^{-1}$.
Hence, by the Chebotarev density theorem, each $\phi_{\sigma}$ is an 
inner automorphism, and our statement follows from the cohomological
Hasse principle for adjoint groups. 

\appendix
\section{On Weil's descent theorem}
Let $k$ be a field, which we will assume for simplicity 
to be of characteristic zero, let $K\supset k$ be an algebraically closed 
field, and let $Y$ be a scheme over $K$.

\begin{Def}
A {\em descent datum} for $Y$ from $K$ to $k$ is a family of isomorphisms 
$\varphi_{\sigma}:{}^{\sigma}Y\isom Y$ for all $\sigma\in\Aut(K/k)$,
which satisfy the relation
\begin{equation} \label{E:cocycle}
\varphi_{\sigma\tau}=\varphi_{\sigma}\circ{}^{\sigma}\varphi_{\tau} \text{ for every } \sigma,\tau\in\Aut(K/k).
\end{equation}
\end{Def}

\begin{Def}

a) A descent datum $\{\varphi_{\sigma}\}_{\sigma\in \Aut(K/k)}$ for $Y$ defines 
natural actions of the group $\Aut(K/k)$ on the underlying topological
space of $Y$ and on the set of rational functions $K(Y)$ by the formulas
$\sigma(y):=\varphi_{\sigma}({}^{\sigma}y)$ and 
$\sigma(f):={}^{\sigma}(\varphi_{\sigma^{-1}}^*(f))$ respectively.

b) Given a descent datum for $Y$ we will say that a subset of $Y$ or a rational function on
$Y$ is {\em defined} or {\em rational} over a subfield $L\subset K$ containing $k$, 
if it is invariant under the action of the subgroup 
$\Aut(K/L)\subset\Aut(K/k)$.  
\end{Def}

\begin{Def}
We will call the descent datum $\{\varphi_{\sigma}\}_{\sigma\in \Aut(K/k)}$ for 
$Y$ {\em continuous} if each rational function on $Y$ is defined over a field, finitely generated
over $k$, or, equivalently, if the stabilizer of each rational
function on $Y$ is open in the Krull topology of $\Aut(K/k)$.
\end{Def}

\begin{Rem}
Every scheme $X/k$ defines a scheme $Y:=X\times_k K$ over
$K$  equipped with a natural continuous descent datum. Conversely, 
we will call a continuous descent datum  
$\{\varphi_{\sigma}\}_{\sigma\in \Aut(K/k)}$ for $Y$ 
{\em effective} if it comes from a certain scheme $X$ over $k$.
In this case, $X$ is unique.
\end{Rem}

The following theorem is essentially due to Weil \cite{We}. 
We decided to include a complete proof of this important fact, 
since it plays a crucial role
for the construction of the canonical models of Shimura varieties 
(see \rco{can}), and since the author is not aware of any written
proof, which uses modern language of algebraic geometry. 

\begin{Rem} \label{R:Milne}
After this work was basically completed we have learned about the 
work of Milne \cite{Mi3}, where he derives similar statement from Weil's theorem.
\end{Rem}
 
\begin{Thm} \label{T:Weil}
Suppose that $K$ has an infinite transcendental degree over $k$ and 
that $Y$ is quasi-projective over $K$. Then every continuous descent 
datum $K/k$ on $Y$ is effective, and the corresponding scheme
$X$ is quasi-projective. Moreover, $X$ is affine or quasi-affine 
whenever $Y$ is so. 
\end{Thm}

\begin{pf}
First we establish some particular cases of the theorem.

\begin{Lem} \label{L:aff}
Theorem holds, if $Y$ is an affine scheme of finite type over $k$.
\end{Lem}
\begin{pf}

a) Let $R$ be the ring of regular functions on $Y$, and let $R_0\subset R$
be the ring of regular functions invariant under the full group 
$\Aut(K/k)$. We have to show that $R_0$  spans $R$ over $K$ and that
it is finitely generated as an algebra over $k$. Since $R$ is finitely 
generated over $K$, the former statement implies the latter.
First we will show that $R$ is spanned by $R'_0:=R^{\Aut(K/\overline{k})}$.

b) Since the algebra $R$ is finitely generated, the continuity assumption 
implies that there exists a subfield $L\subset K$ finitely generated
over $\overline{k}$ such that the algebra $R_1:=R^{\Aut(K/L)}$ spans
$R$ over $K$. Then for every  $\sigma\in\Aut(K/\overline{k})$, 
the conjugate $\sigma(R_1)$ of $R_1$ is the set of all regular
functions on $Y$, defined over $\sigma(L)$. Hence 
for every  $\sigma,\tau\in\Aut(K/\overline{k})$, the algebras
generated by $\sigma(R_1)$ and $\tau(R_1)$ over the composite field 
$Q(\sigma(L)\cdot\tau(L))$ of $\sigma(L)$ and $\tau(L)$ coincide, 
since they both equal to the subalgebra of all regular functions on
$Y$, defined over  $Q(\sigma(L)\cdot\tau(L))$.

c) Suppose now that $\sigma(L)$ and $\tau(L)$ are 
linearly disjoint over $\overline{k}$. Then we claim that the composite 
rings $\tau(R_1)\cdot\sigma(L)$ and $\sigma(R_1)\cdot \tau(L)$ are equal.
\begin{pf}

Let us show, for example, that  $\tau(R_1)\subset\sigma(R_1)\cdot \tau(L)$.
Choose $\mu\in\Aut(K/\overline{k})$ such that $\mu(L)$ is linearly disjoint from $\sigma(L)\cdot\tau(L)$ over
$\overline{k}$ and take $r\in \tau(R_1)$. By b), there exist non-zero elements 
$f\in\tau(L)\cdot\mu(L)$ and $g\in\sigma(L)\cdot\mu(L)$  
such that $rf\in\mu(R_1)\cdot\tau(L)$ and $rfg\in\sigma(R_1)\cdot\tau(L)\cdot\mu(L)$.

Now the statement follows formally from the fact that $\mu(L)$ is linearly disjoint from
$\sigma(R_1)\cdot\tau(L)$ over $\overline{k}$. Namely,
choose a finitely generated subalgebra $A\subset\mu(L)$ over $\overline{k}$ such that 
$f\in\tau(L)\cdot A$, $g\in\sigma(L)\cdot A$ and $rfg\in\sigma(R_1)\cdot\tau(L)\cdot A$. Let $\C{M}$
be a maximal ideal of $A$ such that the reductions $\bar{f}\in\tau(L)$ and $\bar{g}\in\sigma(L)$
of $f$ and $g$ modulo $\C{M}$ are non-zero. Since $r$ belongs to 
$\tau(R_1)$, we have $\overline{r}=r$. Hence
$r\bar{f}\bar{g}\in\sigma(R_1)\cdot\tau(L)$, and the statement
follows, since $\bar{f}$ and $\bar{g}$ are invertible.
\end{pf}

d) Fix now $\sigma\in\Aut(K/\overline{k})$ such that  $\sigma(L)$ is linearly disjoint from $L$ 
over $\overline{k}$, and set $R_2:=R_1\cap\sigma(R_1)$. Choose any $\tau\in\Aut(K/\overline{k})$
such that $\tau(L)$ is linearly disjoint from $L\cdot\sigma(L)$ over
$\overline{k}$. Then $\tau(L)$ is linearly disjoint from $R_1\cdot\sigma(R_1)$ over
$\overline{k}$. Hence c) implies that
$$R_2\cdot\tau(L)=(R_1\cdot\tau(L))\cap(\sigma(R_1)\cdot\tau(L))=
(\tau(R_1)\cdot L)\cap(\tau(R_1)\cdot\sigma(L))=\tau(R_1).$$
In particular, elements of $R_2$ are fixed by any subgroup
$\Aut(K/\tau(L))$, therefore by all of $\Aut(K/\overline{k})$. This means
that $R_2=R'_0$ and $R'_0\cdot L=R_1$. Hence 
$R'_0\cdot K=R_1\cdot K=R$, as claimed. 

e) It remains to show that $R_0=(R'_0)^{\Aut(\overline{k}/k)}$ generates
$R'_0$ as a vector space over $R$, but this is proven in \cite[AG, 14.2]{Bo}.
\end{pf}

\begin{Cor} \label{C:qaff}
 Theorem holds, if  $Y$ is a quasi-affine scheme of finite type over $k$.
\end{Cor}
\begin{pf}
Let $R$ be the ring of regular functions on $Y$. Then the natural
map $Y\to\Spec R$ is an open embedding. But the descent datum on $Y$
defines descent data on affine schemes $\Spec R$ and 
$(\Spec R\sm Y)_{red}$. Therefore the effectivity of the descent for 
quasi-affine schemes is a formal consequence of that for the affine ones.
\end{pf}

\begin{Lem} \label{L:qproj}
If a scheme $Y$ defined over a field $k$ of characteristic zero is 
quasi-projective over its overfield $K$, then it is quasi-projective 
over $k$ itself. 
\end{Lem}

\begin{pf}
Since $Y$ is quasi-projective over $K$, it is quasi-projective
over a certain finitely generated over $k$ subfield $L\subset K$
 (see Step 3 below).
Let $A$ be a finitely generated algebra $A$ over $k$ with fraction field
$L$, then there exists a quasi-projective scheme $\wt{Y}$ over 
$\Spec A$, whose generic fiber is $Y_L$. By \rl{mor} below, the restriction
of  $\wt{Y}$ to a certain open dense subset $U\subset\Spec A$ is 
isomorphic to $Y_U:=Y\times_k U$, so that $Y_U$ is quasi-projective
over $U$. 

By restricting $Y_U$ to a fiber over some closed point of  
$U$, we see that $Y$ is quasi-projective over some finite 
extension $k'$ of $k$, which we may assume is Galois.
This means that $Y_{k'}$ possesses a very ample line bundle $\C{L}$.
Taking the product of the Galois conjugates of $\C{L}$ over $k$, we find a very 
ample line bundle on $Y$, showing therefore that $Y$ is 
quasi-projective over $k$. 
\end{pf}

\begin{Lem} \label{L:mor}
Let $X$ and $Y$ be two schemes of finite type over an irreducible scheme
$Z$. If $X$ and $Y$ have isomorphic generic fibers over $Z$, then they
have isomorphic restrictions to the inverse images of
a certain open dense subscheme of $Z$.
\end{Lem}
\begin{proof}
The assumption implies that $X$ and $Y$ are birationally isomorphic,
therefore they contain open dense subsets 
$U_1\subset X$ and $U_2\subset Y$, isomorphic over $Z$ and 
containing the generic fiber of $X$ and $Y$
respectively. Since $X$ and $Y$ are of finite type over $Z$, 
the projections of $X\sm U_1$ and $Y\sm U_2$ to $Z$ are not Zariski
dense. Therefore there exists an open dense subset 
$U\subset Z$, whose inverse images in 
$X$ and $Y$ are contained in $U_1$ and $U_2$ respectively. Then $U$ is
the required subscheme of $Z$.
\end{proof}

Now we start the proof in the general case.

\begin{Not}
For each rational function $f$ on a scheme $Z$, let $Z^f$ be the
largest open subscheme of $Z$ on which $f$ is regular.  
\end{Not}

{\bf Step 1}. It follows from the above lemmas that it remains to show that each 
point $y\in Y$ has an open affine neighborhood, invariant under the 
action of $\Aut(K/k)$. Moreover, since finite Galois descent is 
effective for quasi-projective schemes (see \cite[Ch. V, 20,
Cor. 2]{Se}), we may replace $k$ by its finite extension.

Since a noetherian scheme is affine if and only if 
the reduced scheme, associated to it is affine 
(see \cite[Ch. III, Ex. 3.1]{Ha}), it will suffice to
assume that $Y$ is reduced.  

As $Y$ has finitely many irreducible components, the continuity
of descent datum imply that each irreducible component of $Y$ 
is defined over a finite extension of $k$. Therefore
we can replace $Y$ by its connected component and $k$ by its 
finite extension, that is, we can  assume that $Y$ is integral.

 Let $Y_0$ is the set of all points of $Y$ which does not have a 
required neighborhood. Then $Y_0$ is a closed subset of $Y$, rational
over $k$, and we want to show that it is empty. 
If not, let $y_0$ be the generic point of some irreducible component of 
$Y_0$. As above, $y_0$ is rational over a finite extension of $k$,
hence we may assume after enlarging $k$ that $y_0$ is rational over $k$.

To get a contradiction, we will find a rational function $f$ on $Y$,
defined over $\overline{k}$, for which $Y^f$
is an affine neighborhood of $y_0$. The continuity assumption will
then imply that $f$ and, therefore, $Y^f$ are defined 
over a finite Galois extension of $k$, therefore the intersection
of all Galois conjugates of $Y^f$ over $k$ will give the required affine
neighborhood of $y_0$, contradicting to our choice of $y_0$.

{\bf Step 2}. First we claim that
 field $M:=K(Y)^{\Aut(K/\overline{k})}$ generates $K(Y)$ as a
 field over $K$.
\begin{pf}
 The argument is very similar to that of \rl{aff}.
Since the field $K(Y)$ is finitely generated over $K$, the continuity
of the descent datum implies that there exists a subfield $L\subset K$ 
finitely generated  over $\overline{k}$ such that the field $M_L$ of 
rational functions on $Y$ which are  defined over $L$ generates $K(Y)$ as a 
field over $K$.

Choose now $\sigma,\tau\in\Aut(K/\overline{k})$ as in the part c) 
of the proof of \rl{aff}, and set $M':=M_L\cap\sigma(M_L)$.
Then we claim that the intersection

\begin{equation} \label{E:int}
Q(M_L\cdot\tau(L))\cap Q(\sigma(M_L)\cdot\tau(L))\subset 
Q(M_L\cdot\sigma(M_L)\cdot\tau(L)) 
\end{equation}
equals $Q(M'\cdot\tau(L))$. In fact, let $f$ belongs to the intersection. 
Choose now a finitely generated over $\overline{k}$ subalgebra 
$A\subset\tau(L)$, whose fraction field is $\tau(L)$. Then for almost all maximal ideals $\C{M}\subset A$ the
reduction $\bar{f}$ of $f$ modulo $\C{M}$ is defined. Since every such an $\bar{f}$ necessary belongs to
$M_L\cap\sigma(M_L)=M'$, we get that $f\in Q(M'\cdot\tau(L))$, as was claimed.

As in the part b) of the proof of \rl{aff}, the intersection (\ref{E:int}) is equal to 
$Q(\tau(M_L)\cdot L)\cap Q(\tau(M_L)\cdot\sigma(L))=\tau(M_L)$. The final arguments 
of part d) of the proof of \rl{aff}  can now be applied to our situation without any changes.
\end{pf}

{\bf Step 3}. Next we show that the descent datum 
$\{\varphi_{\sigma}\}_{\sigma\in\Aut(K/L)}$ on $Y$ is effective for some
finitely generated over $\overline{k}$ subfield $L\subset K$.

Indeed, choose an embedding of $Y$ into a certain projective 
space $\B{P}_K^n$ with homogeneous coordinates $(x_0:...:x_n)$.
After a possible change of coordinates we may assume that all the
rational functions $f_{ij}:=x_i/x_j$ on $\B{P}^n_K$ define rational 
functions on $Y$.

Let $f_1,...,f_l$ and $g_1,...,g_k$ be homogeneous polynomials
in $n+1$ variables $x_0,...,x_n$ such that 
\begin{equation} \label{E:eq}
Y=\{\bar{x}\in\B{P}_K^n\,|\, f_1(\bar{x})=...=f_l(\bar{x})=0\}\sm
\{\bar{x}\in\B{P}_K^n\,|\,g_1(\bar{x})=...=g_k(\bar{x})=0\},
\end{equation}
and let $L$ be a finitely generated extension of $\overline{k}$,
containing all the coefficients of the $f_i$'s, that of the $g_j$'s and
fields of rationality of the $f_{ij}$'s. Then equations (\ref{E:eq})
define the required quasi-projective variety $Y'$ over $L$.

 
{\bf Step 4}.
To find the required function $f$ from Step $1$ we first choose
a   non-zero rational  function $g$  on $Y$ of the form 
$h_2^{\text{deg}\,h_1}/h_1$, where $h_1$ is one of the 
 $g_i$'s from Step $3$, which does not vanish at $y_0\in Y$, and  
$h_2\in L[x_0,...,x_n]$ is a linear polynomial such that 
the closed subset $\{z\in Y\,|\,h_1(z)=h_2(z)=0\}$ has a codimension 
two in $Y$. Then $g$ is defined over $L$, and  
$Y^g=\{z\in Y|h_1(z)\neq 0\}$ is an affine neiborhood of $y_0$.

Take any $\sigma\in\Aut(K/k)$ such that $\sigma(L)$ is linearly 
disjoint from $L$ over $\overline{k}$, and put 
$h:=\sigma(g)\subset K(Y)$.
Then $h$ is a rational function on $Y$, defined over $\sigma(L)$,
and $Y^h$ is also an affine neighborhood
of $y_0$. The result from Step $2$ implies that $h$ belongs to $Q(M\cdot\sigma(L))$.

 Note that since $y_0\in Y$ is rational over $k$, it descends to a point 
$y'_0\in Y'$. Choose now a finitely generated over $\overline{k}$ subalgebra 
$A\subset\sigma(L)$ whose quotient field is $\sigma(L)$.
Then $h$ can be considered as a rational function on 
$Y'\times_{\overline{k}} \Spec\,A$ regular at $y'_0\times\Spec\,\sigma(L)$.
There exists an open affine dense subset $V\subset\Spec\,A$ such that 
for all closed points $x$ of $V$ the restriction $h_x$ of $h$ to 
$Y'\times\{x\}\cong Y'$ gives a rational function on $Y'$, 
regular at $y'_0$. Since each $h_x$ 
will automatically belong to $M$, it remains to find an $x$ with 
$(Y')^{h_x}$ affine.

{\bf Step 5}. Since $Y^h$ is affine, we can shrink $V$ so that
$U:=(Y'\times_{\overline{k}}V)^h$ will be affine as well (use
\rl{mor}). Hence it will suffice to find an $x$ such that 
$U_x:= U\cap (Y'\times\{x\})\subset Y'$ equals $(Y')^{h_x}$.
Denote $(Y'\times_{\overline{k}}V)^{h^{-1}}$ by $U'$, and for each 
closed point $x$ of $V$ put  $U'_x:= U'\cap (Y'\times\{x\})\subset Y'$.
Then we have the following inclusions:

$$U'_x\sm U_x\subset Y'\sm (Y')^{h_x}\subset Y'\sm U_x.$$
Indeed, the last inclusion is obvious, and the first follows from 
the observation that function $h^{-1}$  is regular and vanishes at all 
points of $U'_x\sm U_x$, therefore $h_x$ can not be regular at such 
points.

Our choice of $g$ implies that $Y\sm(Y^g\cup Y^{g^{-1}})$ and, 
therefore, $Y\sm(Y^h\cup Y^{h^{-1}})$ are of codimension two in $Y$.
 Since the dimension of the fibers is generically constant 
(see \cite[I, $\S$6.3, Thm. 7]{Sh}), there exists a closed point 
$x'\in V$ such that $Y'\sm (U_{x'}\cup U'_{x'})$ has a codimension two 
in $Y'$.
Since $U_{x'}$ is affine, \rl{codim} below implies that 
the set $U'_{x'}\sm U_{x'}=[Y'\sm U_{x'}]\sm [Y'\sm (U_{x'}\cup U'_{x'})]$
is Zariski dense in $Y'\sm U_{x'}$.
Since $Y'\sm (Y')^{h_x'}$ is a closed subset of $Y'$, we get 
$ (Y')^{h_{x'}}= U_{x'}$, thus completing the proof of \rt{Weil}.
\end{pf}

\begin{Lem} \label{L:codim}
If $U$ be an open and dense affine subset of a separated irreducible
noetherian scheme $X$, then the complement of $U$ in $X$ has a pure 
codimension one (that is, all the irreducible components of $X\sm U$ 
have codimension one in $X$).
\end{Lem}
\begin{pf} (compare \cite[II, Prop. 3.1]{Ha2}, where a cohomological 
proof of this fact is given).
First we notice that we can replace
$X$ by the associated reduced scheme, thus assuming that
$X$ is reduced (use \cite[Ch III, Ex. 3.1]{Ha}). Since 
the normalization morphism is finite (hence
affine), we may 
replace $X$ and $U$ by their normalizations, and thus to assume 
that $X$ is normal. Then  
\cite[Thm. 38]{Mat} implies that for every 
rational function $f$ on $X$ the complement of $X^f$ has a pure 
codimension one in $X$. Since $U$ is affine, it equals 
the intersection of the $X^f$'s taken over the set of all regular
functions $f$ on $U$, considered as rational functions on $X$. 
Hence the complement of $U$ in $X$ is the union
of the compliments of the $X^f$'s. Therefore it 
has also a pure codimension one in $X$, as was claimed.
\end{pf}

The following proposition gives a convenient way for checking that
the descent datum is continuous.

\begin{Prop} \label{P:cont}
(General principle) Let $K$ and $k$ be as in the Theorem, and let
$Y$ be a scheme of finite type over $K$, equipped with a descent
datum $\{\varphi_{\sigma}\}_{\sigma\in\Aut(K/k)}$.

Suppose that $Y$ can be equipped with some extra datum $\C{D}$ of
finite type over $Y$ such that no non-trivial automorphism of $Y$ 
extends to an automorphism of $\C{D}$.

If  $\{\varphi_{\sigma}\}_{\sigma}$ can be extended to a descent datum of 
$\C{D}$, then $\{\varphi_{\sigma}\}_{\sigma}$ is automatically continuous.
\end{Prop}

\begin{pf}
Since $\C{D}$ is of finite type, it has a structure over some 
finitely generated  over $k$ subfield $L\subset K$.
This structure define a continuous descent datum 
$\{\psi_{\tau}:{}^{\tau}\C{D}\isom\C{D}\}_{\tau}$ from $K$ to $L$.

Let $\{\wt{\varphi}_{\sigma}\}_{\sigma}$ be an extension
of $\{\varphi_{\sigma}\}_{\sigma}$ to the descent datum
of $\C{D}$. Then for each $\tau\in\Aut(K/L)$ the composit map
 $\wt{\varphi}_{\tau}^{-1}\circ\psi_{\tau}$
is an automorphism of $\C{D}$, hence, by our assumption, 
it induces the identity automorphism on $Y$. In particular,
 the isomorpism 
${}^{\tau}Y\isom Y$ induced by  $\psi_{\tau}$ coincides with
$\varphi_{\tau}$ for all $\tau\in\Aut(K/L)$. Since the descent datum 
$\{\psi_{\tau}\}_{\tau\in\Aut(K/L)}$ was chosen to be continuous, and 
since $L$ is finitely generated over $k$, our original descent datum
$\{\varphi_{\sigma}\}_{\sigma}$ is continuous, as was claimed.
\end{pf}

\begin{Cor} \label{C:aut}
 Let $K$ and $k$ be as in the Theorem. Suppose that we are given
a data consisting of a family $\{X_{\al}\}_{\al\in\C{A}}$ of 
quasi-projective varieties over $K$ and a family 
$\{f_{\beta}\}_{\beta\in\C{B}}$ of morphisms between them.

Suppose that the automorphism group of each $X_{\al}$ is finite
and that there is no non-trivial families of automorphisms
 $\{\psi_{\al}\in \Aut(X_{\al})\}_{\al\in\C{A}}$ commuting with 
all the $f_{\beta}$'s.

Then every descent datum of the whole family 
$\{X_{\al}\}_{\al\in\C{A}}$, which commutes with all 
the $f_{\beta}$'s, is effective. 
\end{Cor}
\begin{pf}
For every finite subset $I\subset\C{A}$ let $J(I)\subset\C{B}$ 
be the set of all the $\beta$'s in $\C{B}$ such that $f_{\beta}$
is a map from  $X_{\al_1}$ to $X_{\al_2}$ with $\al_1, \al_2\in I$.
Denote now by $\Aut_I$ the group of families 
$\{\psi_{\al}\in \Aut(X_{\al})\}_{\al\in I}$, 
commuting with  all the $f_{\beta}$'s with $\beta\in J(I)$.

For each finite sets $I_1\subset I_2$ we have a natural
group homomorphism $\Aut_{I_2}\to\Aut_{I_1}$. 
In this way $\{\Aut_I\}_I$ form a projective system, and our assumption 
says that the inverse limit of the $\Aut_I$'s is the trivial group.

Fix now an $\al\in\C{A}$. Since $\Aut(X_{\al})$ is finite, 
the triviality of the inverse limit implies that there exists a finite 
subset $I_{\al}$ of $\C{A}$ containing $\al$ such the natural homomorphism
 $\Aut_{I_{\al}}\to\Aut_{\{\al\}}$ maps all the group to the 
identity. Since $I_{\al}$ is finite and since 
$\Aut(X_{\gm})$ is finite for each $\gm\in I_{\al}$, there exists 
therefore a finite subset $J_{\al}\subset J(I_{\al})$ such that
no non-trivial automorphism of $X_{\al}$ extends
to an automorphism of a datum  consisting of 
quasi-projective varieties $\{X_{\gm}\}_{\gm\in I_{\al}}$ and morphisms
$\{f_{\beta}\}_{\beta\in J_{\al}}$. \rp{cont} and \rt{Weil}
now imply that the corresponding descent datum on $X_{\al}$
is effective. Since $\al$ was chosen arbitrary, we get the assertion.
\end{pf}

\end{document}